\documentclass[a4paper,11pt]{amsart}

\usepackage[textwidth=15cm, textheight=21cm]{geometry}

\usepackage{amsmath,amssymb,amsthm,amscd,amsfonts}
\usepackage[hyphenbreaks]{breakurl}
\usepackage{mathrsfs}
\usepackage{latexsym}
\usepackage[all]{xy}
\usepackage{verbatim}
\usepackage[usenames, dvipsnames]{color}
\usepackage{multirow}
\usepackage{url}
\usepackage{mathdots}
\usepackage{MnSymbol}
\usepackage{stmaryrd}
\allowdisplaybreaks

\usepackage{stmaryrd}
\usepackage{tikz}
\usetikzlibrary{%
    matrix,%
    calc,%
    arrows,%
    automata,%
    chains,%
    positioning,%
    scopes
}
\tikzset{>=stealth',every on chain/.append style={join},
         every join/.style={->}}

\newcommand{\field}[1]{\mathbb{{#1}}}
\newcommand{\N}{\field{N}}
\newcommand{\Z}{\field{Z}}
\newcommand{\Q}{\field{Q}}

\newcommand{\C}{\field{C}}

\newcommand{\F}{\field{F}}
\newcommand{\Fq}{\field{F}_q}
\newcommand{\Fi}{\field{F}_\infty}
\newcommand{\Falg}{\overline{\field{F}}_q}

\newcommand{\Zp}{\field{Z}_p}
\newcommand{\Qp}{\field{Q}_p}

\newcommand{\Si}{\field{S}_{\infty}}
\newcommand{\Sv}{\field{S}_{\nu}}
\newcommand{\If}{\field{I}_{F}}
\newcommand{\Gal}{{\rm Gal}}
\newcommand{\Ker}{{\rm Ker}}

\newtheorem{teo}{Theorem}[section]              
\newtheorem*{teo2}{Theorem}                        
\newtheorem{prop}[teo]{Proposition}         

\newtheorem{lemma}[teo]{Lemma}

\newtheorem{defi}[teo]{Definition}
\newtheorem{osse}[teo]{Remark}

\title{Stickelberger series and Main Conjecture for function fields}

\author{Andrea Bandini}

\address{{\sc Andrea Bandini}: Universit\`a degli Studi di Pisa \\
	Dipartimento di Matematica\\
	Largo Bruno Pontecorvo, 5 \\
	56127 Pisa - Italy
}

\email{andrea.bandini@unipi.it}

\author{Edoardo Coscelli}

\address{{\sc Edoardo Coscelli}: 
}

\email{coscelliedoardo@gmail.com}

\begin{document}

	\begin{abstract}
	\noindent	Let $F$ be a global function field of characteristic $p$ with ring of integers $A$ and let $\Phi$ be a Hayes module on the Hilbert 
	class field 	$H_A$ of $F$. 
		We prove an Iwasawa Main Conjecture for the $\Z_p^\infty$-extension $\mathcal{F}/F$ generated by the $\mathfrak{p}$-power
		torsion of $\Phi$ ($\mathfrak{p}$ a prime of $A$). The main tool is a Stickelberger series whose specialization
		provides a generator for the Fitting ideal of the class group of $\mathcal{F}$. Moreover we prove that the same series, evaluated at complex 
		or $\mathfrak{p}$-adic	characters, interpolates the Goss Zeta-function or some $\mathfrak{p}$-adic $L$-function, thus providing the link 
		between the algebraic structure (class groups) and the analytic functions, which is the crucial part of Iwasawa Main Conjecture.
	\end{abstract}


	
\maketitle
	
	\section{Introduction}

Arithmetic properties of motives defined over a global field are a central theme of modern number theory
and one of the main topic is their relation with (or interpretation as) special values of $L$-functions.
Iwasawa theory combines an algebraic approach, which studies the variation of motives in $p$-adic towers
as modules over an algebra of power series, with the definition of analytical $p$-adic $L$-functions,
which interpolate special values of more classical $L$-functions, thus providing a link between the 
two aspects of the theory. The Main Conjecture (IMC) predicts (in a growing number of incarnations) 
a deep relation between $p$-adic $L$-functions and a generator of the characteristic or Fitting ideal
of the algebraic structure we deal with (class groups, Selmer groups and so on).

\noindent This topic has been extensively studied since the first proof of a Main Conjecture by Mazur and Wiles (in \cite{Mazur Wiles}) for the class groups of the $\Z_p$-cyclotomic extension of $\Q$, 
but mostly for $\Z_p$ or $\Z_p^d$-extensions of number fields (see e.g. \cite{KatoICM}).
Iwasawa theory for function fields in positive characteristic is relatively new and, even if
some instances of the IMC for $\Z_p^d$-extensions were proved by Crew in \cite{Crew}  and 
Burns, Lai and Tan in \cite{Bur}, we believe that the true analog of the Mazur-Wiles theorem
is the one recently presented in \cite{ABBL}. We briefly explain the setting of that result:
let $F=\F_q(t)$ be the rational function field which plays the role of $\Q$, fix the prime 
$\frac{1}{t}=\infty$ and consider $A=\F_q[t]$ as the ring of integers of $F$, i.e. the functions 
regular outside $\infty$. Let $\mathfrak{p}$ be a prime of $A$ and let $F_n$ be the field generated
over $F$ by the $\mathfrak{p}^{n+1}$-torsion of the Carlitz module associated to $A$ 
(see \cite[Section 2]{ABBL} for a quick review): put $\mathcal{F}=\cup F_n$, then 
$\Gal(\mathcal{F}/F)\simeq \Gal(F_0/F)\times \Gal(\mathcal{F}/F_0)\simeq\Delta\times\Z_p^\infty$,
where $\Delta$ has finite order prime with $p$. The field $\mathcal{F}$ is the analog of the 
cyclotomic $\Z_p$-extension of a number field constructed with $p$-power roots of unity and it is
this analogy which led us to believe that we have to work with $\mathcal{F}$ instead of 
$\Z_p^d$-extensions for some finite $d$. 

\noindent This provides one of the main differences with the number field setting where there are
no $\Z_p^\infty$-extensions: here our algebraic structures (namely class groups) will be modules over 
the non-noetherian algebra $\Lambda(\mathcal{F})=W\llbracket \Gal(\mathcal{F}/F_0)\rrbracket$ 
($W$ an appropriate finite extension of $\Z_p$) and,
having no structure theorem for them, we shall use a limit process to describe their Fitting ideals.  

\noindent Let $\mathcal{C}\ell^0_n$ be the $p$-part of the group of divisor classes of degree 0
of $F_n$ and, for any complex character $\chi$ of $\Delta$, let $\mathcal{C}\ell^0_n(\chi)$ be its
$\chi$-part, i.e. the one on which $\Delta$ acts via $\chi$. Put $\mathcal{C}\ell^0_\infty(\chi)$
for the inverse limit with respect to norms of the $\mathcal{C}\ell^0_n(\chi)$, then 
\cite[Theorems 5.1 and 5.2]{ABBL} prove that for any nontrivial $\chi$ the 
$\Lambda(\mathcal{F})$-module $\mathcal{C}\ell^0_\infty(\chi)$ is finitely generated and torsion
and its Fitting ideal is generated by a specialization of the $\chi$-part of a Stickelberger series
$\Theta(X)\in \Z[\Gal(\mathcal{F}/F)]\llbracket X\rrbracket$.

\noindent The same Stickelberger series is also used to interpolate (via complex or 
$\mathfrak{p}$-adic characters) the Artin $L$-function, the Goss Zeta-function and a $\mathfrak{p}$-adic
$L$-function (\cite[Section 3]{ABBL}), thus providing the link between the algebraic structure
of $\mathcal{C}\ell^0_\infty(\chi)$ and various $L$-functions, i.e. the IMC.

\subsection{Our results} Our goal here is to extend the main results of \cite{ABBL}
to a general function field $F$ of characteristic $p$, i.e. a field of transcendence degree 1 over
some finite field $\F_q$ or, equivalently, the function field of a smooth projective curve defined
over $\F_q$. Fix a prime $\infty$, let $A$ be the ring of integers for $F$ and write $H_A$ for the\
Hilbert class field of $A$, i.e. the maximal abelian unramified extension of $F$ where $\infty$
is totally split. 

\noindent To deal with this setting we need to consider extensions generated by the $\mathfrak{p}$-torsion
of a sign-normalized rank 1 Drinfeld module (or Hayes module, see Definition \ref{DefHayesModule})
and decompose Iwasawa modules in eigenspaces with respect to characters whose action on the inertia groups of $\infty$
is crucial. We replace the odd and even characters of the case $F=\F_q(t)$ (where $H_A=F$) with three types (Definition 
\ref{DefCharType}): we shall provide a complete study for the first two types, the characters of type 3 need a
different treatment (see Remark \ref{RemType3}).

\noindent Moreover the presence of a nontrivial Hilbert class field $H_A$ enters in the definition of the Goss Zeta-function
(in particular in the exponentiation of ideals) and we have to take it into account in all our interpolation formulas and 
in the definition of our $\mathfrak{p}$-adic $L$-series.

\noindent Fix a finite set $S$ of primes of $F$ and let $F_S$ be the maximal abelian extension of 
$F$ unramified outside $S$ with $G_S:=\Gal(F_S/F)$. For any prime $\nu$, let $d_\nu$ be its degree
and, if $\nu\not\in S$, denote by $\phi_\nu$ the Frobenius of $\nu$ in $G_S$. The leading role in this paper 
will be played by the Stickelberger series
\[ \Theta_S (X) = \prod_{\nu \not\in S}\left(1 - \phi_{\nu}^{-1} X^{d_\nu}\right)^{-1}
\in \Z [G_S] \llbracket X \rrbracket \,.\] 

\subsubsection{Analytic side} In Sections \ref{SecStickSer}, \ref{SecStickGossZeta} and 
\ref{SecSticknuadicL} we prove convergence properties for $\Theta_S(X)$ and test it against
characters $\Psi:\Z[G_S]\longrightarrow L$ to provide interpolation formulas for\begin{itemize}
	\item the Artin $L$-function $L_F$ (for $L=\C$, Theorem \ref{ThmCollegamentoArtinL});
	\item the Goss Zeta-function $\zeta_A$ (for $L=\C_\infty$, Theorem 
	\ref{teorema stickelberger-goss});
	\item various $\nu$-adic Zeta-functions $\zeta_\nu$ (for $L=\C_\nu$,
	Theorem \ref{teorema interpolazione v-adica}),
	\end{itemize}
where $\C_\eta$ is the completion of an algebraic closure of the completion $F_\eta$ of $F$
at a prime $\eta$ (for both cases $\eta$ finite or infinite).  Then we construct a $\nu$-adic $L$-series which, via its relation with $\Theta_S(X)$,
interpolates the values of $\zeta_\nu$ at integers for $\nu\in S-\{\infty\}$ (and a bit more, see
Theorem \ref{teorema interpolazione v-adica2}).

As an example we mention here the interpolation we obtain for the Goss Zeta-function 
\[ \zeta_A(s) = \sum \mathfrak{a}^{-s} \qquad 
\text{ with } s \in \Si:=\C_\infty^\times \times \Zp \, , \]
where the sum is taken over all the nonzero ideals $\mathfrak{a}$ of $A$. 
This function represents the natural analogue of the Dedekind Zeta-function of a number field.
For any $y\in \Z_p$ we define a continuous character $\Psi_y : G_S \rightarrow \C_\infty^\times$ 
via
\[\Psi_y (\phi_\nu) = N \left( \langle \nu \rangle_\infty^{-1} \right)^{y/f} \]
(details in Section \ref{SecInterpolation of the Goss Zeta-function}, here it suffices to say
that $N$ is a norm, $f$ is a power of $p$ and $\langle \nu \rangle_\infty$ is a 1-unit so that taking the $\frac{y}{f}$-th power makes sense). 
For every $s=(x,y) \in \Si$ we have
	\begin{equation}\label{EqIntro}
	\Psi_y \left( \Theta_S (X) \right) (x) = \zeta_A (-s) \prod_{\nu \in S \,,\, \nu \neq \infty}(1-\nu^s) \,.
	\end{equation}

\subsubsection{Algebraic side} 
Let $\Phi : A \rightarrow H_A \{ \tau \}$ be a Hayes module (more details on all objects mentioned here are in 
Section \ref{SecStickClassGroups}): fix a prime $\mathfrak{p}$ of $F$ with degree $d_\mathfrak{p}$ and denote 
by $F_n$ the extension of $H_A$ generated by the $\mathfrak{p}^{n+1}$-torsion of $\Phi$. The field $F_n$ is an abelian
Galois extension of $F$ ramified only at $\mathfrak{p}$ and $\infty$. These fields form an 
Iwasawa tower
\[
F \subset H_A \subset F_0  \subset F_1 \subset \dots \subset F_n \subset \dots \subset \bigcup_{n \in \N}F_n =: \mathcal{F}
\]
and, if we put $\Gamma_n = \Gal(F_n / F_0)$, we have 
\[
\Gamma_\infty := \Gal(\mathcal{F} / F_0) = \lim_\leftarrow \Gal(F_n / F_0) \simeq \Zp^{\infty}.
\]
Let $C_n:=\mathcal{C}\ell^0(F_n)$ be the $p$-part of the class group of degree zero divisors 
of $F_n$: it is a $\Zp[\Gamma_n]$-module in a natural way and the $C_n$ form a projective system 
with respect to the norm maps. We define
$\displaystyle{ C_\infty:= \lim_\leftarrow C_n}$, which is a module over the 
non-noetherian \emph{Iwasawa algebra} $\Zp \llbracket \Gamma_\infty \rrbracket$. 
To study the structure of $C_\infty$ we consider $\chi$-parts with respect to 
the characters of the group $G_0:=\Gal(F_0/F)$ which acts naturally on $C_n$ and $C_\infty$. 
We need two simplifying assumptions (which are not required for the analytic part of the theory):
\begin{itemize}
\item $\deg(\infty)=d_\infty=1$, this ensures that all extensions we work with are geometric but it is not restrictive because 
we can always reduce to this case by enlarging the constant field of $F$;
\item the class number of degree zero divisors $h^0(F)$ is prime with $p$, this is necessary to avoid characters of order $p$,
but it is not too restrictive.
\end{itemize}
We extend all scalars to $W=\Zp[\zeta]$, the Witt ring generated by a root of unity $\zeta$ of order $|G_0|$.
Let $S=\{\mathfrak{p},\infty\}$, let $\Theta_\infty(X)$ be the projection of the Stickelberger 
series $\Theta_S(X)$ to $\Z[ \Gamma_\infty \times G_0]\llbracket X \rrbracket$ and write
 $\Theta_\infty(X,\chi) = \chi(\Theta_\infty)(X) \in W[\Gamma_\infty] \llbracket X \rrbracket$
as its $\chi$-part. 
We put
\begin{displaymath}
\Theta_\infty^{\sharp}(X,\chi) =
\left\{ \begin{array}{ll}
\Theta_\infty(X,\chi)                         &   \ \mbox{if } \chi \mbox{ is of
	type } 1 \,, \\
\ & \\                                                                    
\displaystyle{\frac{\Theta_\infty(X,\chi)}{1-X} }   &   \ \mbox{if } \chi \mbox{ is of
	type } 2 \,.
\end{array} \right.
\end{displaymath}
Computing Fitting ideals for the $C_n$ and working our way through a limit process we prove

\begin{teo2}{\rm (Theorems \ref{ThmClassGrFinGen} and \ref{TeoFittIdealClGrIwasawaModule})}
	Let $\chi$ be a character of type $1$ or $2$. Then $C_\infty(\chi)$ is a finitely generated
	module over the Iwasawa algebra $\Lambda:=W\llbracket \Gamma_\infty \rrbracket$, and
	\begin{displaymath}
	{\rm Fitt}_{\Lambda} \left(C_\infty(\chi) \right) =
	\left( \Theta_\infty^{\sharp}(1,\chi)\right) \,.
	\end{displaymath}
\end{teo2}

This is what we call IMC for this setting, the relations with (special values of) 
$\mathfrak{p}$-adic $L$-functions are provided by the analytic interpolation properties
mentioned before.

\begin{osse}\label{RemType3}
For the characters of type 3 we are only able to compute the Fitting ideal of a dual of $C_n$ and
it is often non principal: we have no arithmetic interpretation (from the point of view of Iwasawa 
theory) for this situation yet so we decided to present it in a different paper (see \cite[Section 3]{BBC}).
\end{osse}

\subsection{Setting and notations}

\noindent\begin{enumerate}
	\item[$\bullet$] $F$ is a global function field of characteristic $p>0$, i.e. a finite algebraic extension of a field of
	transcendence degree 1 over a finite field $\F_{p^r}:=\Fq$ which we call the {\em constant field} of $F$.
	A more geometric interpretation would be to consider $F$ as the function field of a smooth projective curve
	$X$ defined over $\Fq$;
	\item[$\bullet$] $\infty$ is a fixed place of $F$ and $A$ is the subring of $F$ of the elements regular outside the place $\infty$;
	\item[$\bullet$] for any place $\nu$ of $F$ (including $\infty$), $F_{\nu}$ is the completion of $F$ at $\nu$. Its ring of integers
	will be denoted by $O_\nu$ and $U_1(\nu)$ will be the group of $1$-units of $F_\nu$.
	The residue field $O_\nu/(\nu):=\F_\nu$ is a finite extension of $\Fq$ of degree $d_\nu := [\F_\nu : \Fq]$
	(also called the {\em degree} of $\nu$), its order will be denoted by $\textbf{N}\nu := q^{d_\nu}$;
	\item[$\bullet$] $v_\nu: F_\nu \rightarrow \Z$ is the (canonical) valuation at $\nu$ and $\pi_\nu$ will denote a fixed uniformizer
	for $F_\nu$, i.e. an element with $v_\nu(\pi_\nu)=1$;
	\item[$\bullet$] the degree of a fractional ideal $\mathfrak{a} = \displaystyle{\prod_{\nu \neq \infty} \nu^{n_\nu}}$ of $A$ is the 
	quantity $\deg(\mathfrak{a}) = \displaystyle{\sum_{\nu \neq \infty} n_\nu d_\nu}$.
\end{enumerate}

\section{Stickelberger series and Artin $L$-function}\label{SecStickSer}
In this first section we introduce the main object of both the analytic and algebraic sides of the theory: the {\em Stickelberger series}. We
provide here its convergence properties and its relation with complex Artin $L$-functions. In the subsequent sections it will appear as
an incarnation of various $L$ or Zeta-functions and, on the algebraic side, as a generator for Fitting ideals of class groups.

\subsection{Stickelberger series}
Let $S$ be a nonempty finite set of places of $F$ that contains $\infty$ and denote by $F_S$ the maximal abelian extension of $F$
unramified outside $S$ and with $G_S := \Gal(F_S / F)$ its Galois group. For every place
$\nu \not\in S$, let $\phi_\nu$ be the Frobenius at $\nu$, i.e. the unique
element of $G_S$ that satisfies
\[ \phi_\nu (x) \equiv x^{\textbf{N}\nu} \pmod{\tilde{\nu}} \]
for every $x \in F_S$, where $\tilde{\nu}$ is any place of $F_S$ lying above $\nu$. 
The extension $F_S / F$ is unramified at $\nu$, so
the decomposition group of $\nu$ in $G_S$ is pro-cyclic and topologically generated by $\phi_\nu$.

\begin{defi}\label{DefStick}
    The Stickelberger series of $S$ is the power series defined by the Euler product
	\[ \Theta_S (X) = \prod_{\nu \not\in S}\left(1 - \phi_{\nu}^{-1} X^{d_\nu}\right)^{-1}
		\in \Z [G_S] \llbracket X \rrbracket \,.\]
\end{defi}

\noindent For an equivalent formula, recall that every Euler factor $e_\nu(X):= 1 - \phi_{\nu}^{-1} X^{d_\nu}$ is invertible in 
$\Z [ G_S ] \llbracket X \rrbracket$ and let $\mathcal{I}_S$ be 
the set of fractional ideals of $A$ with support outside $S$ (recall $\infty \in S$). Let
$\phi_{\mathfrak{a}}$ be the Artin symbol for $\mathfrak{a} \in \mathcal{I}_S$, i.e.
$\displaystyle{\mathfrak{a}=\prod_{\nu\not\in S} \nu^{n_\nu} \mapsto \phi_{\mathfrak{a}} = \prod_{\nu\not\in S} \phi_{\nu}^{n_\nu}}$,
then
\[ \Theta_S (X) = \sum_{\mathfrak{a} \in \mathcal{I}_S ,\, \mathfrak{a} \geqslant 0}
                    \phi_{\mathfrak{a}}^{-1} X^{\deg \mathfrak{a}} =
\sum_{n\geqslant 1} \bigg( \sum_{\substack{\mathfrak{a} \in \mathcal{I}_S ,\, \mathfrak{a} \geqslant 0 \\ \deg \mathfrak{a}=n}}
                    \phi_{\mathfrak{a}}^{-1} \bigg) X^n \]
(where $\mathfrak{a} \geqslant 0$ denotes the integral ideals of $A$).
For any $n$ there exists only a finite number of primes $\nu$ with degree equal to $n$, hence
the series on the right is an element of $\Z [ G_S ] \llbracket X \rrbracket$.

\subsection{Artin $L$-functions}\label{SecArtinL}
Let $K/F$ be a finite subextension of $F_S$ whose Galois group is $G$ and let $S_K \subseteq S$
be the set of ramified places together with $\infty$ (in particular $S_F=\{\infty\}$). 

For every prime $\nu$ let $\phi_{K,\nu} \in G$ be its Artin symbol: if $\nu\not\in S$, then 
$\phi_{K,\nu}$ is the image of $\phi_\nu$ via the canonical projection
$G_S \twoheadrightarrow G$, while for a ramified prime $\nu$, any lifting 
of its Frobenius in $G/I(\nu)$ (where $I(\nu)$ is the inertia at $\nu$ in $G$) will do.  
For every complex character $\chi$ of $G$, i.e. an element of $\text{Hom}(G, \C^{\times})$, we put
$\chi(\nu) = \chi(\phi_{K,\nu})$ 

\begin{defi}\label{DefArtinL}
The {\em Artin $L$-function} associated to $(K, \chi)$ is the complex variable function 
\[     L_{K/F}(s, \chi) = \prod_{\nu\not\in S_K} \left(1 - \chi(\nu) ({\bf N}\nu)^{-s} \right)^{-1},
    \quad \text{for } \mathfrak{Re}(s) > 1  \]
(where the condition $\mathfrak{Re}(s) > 1$ guarantees convergence).
\end{defi}

Our goal in this section is to provide a link between $L_{K/F}(s, \chi)$ and $\Theta_S (X)$. In what follows we present formal
equalities and avoid mentioning the radius of convergence: to be on the safe side one can always assume $\mathfrak{Re}(s)> 1$. 

\noindent Let $\Psi : G_S \rightarrow \C^{\times}$ be a continuous character of $G_S$, i.e. a continuous homomorphism
with respect to the natural topologies: it induces a ring homomorphism $\Z [ G_S ] \llbracket X \rrbracket \rightarrow \C \llbracket X \rrbracket$
(still denoted by $\Psi$).

\begin{teo}\label{ThmCollegamentoArtinL}
    \begin{enumerate}
    \item[(a)] Let $K/F$, $G$ and $\chi$ be as above, then there exists a continuous character $\Psi$ of $G_S$
        such that
        \begin{equation}\label{rappresentazione complessa}
            \Psi \left( \Theta_S (X) \right) (q^{-s}) =
            L_{K/F}(s, \chi^{-1}) \prod_{\nu \in S - S_K}\left(1 - \chi^{-1}(\nu)q^{-s d_{\nu}} \right) \,.
        \end{equation}
    \item[(b)] Let $\Psi$ be a continuous character of $G_S$. Then there exists $K\subset F_S$, finite over $F$ and a complex character $\chi$ of 
       $\Gal(K/F)$ such that equation \eqref{rappresentazione complessa} holds.
    \end{enumerate}
\end{teo}

\begin{proof}
(a) Let $\pi_K$ be the canonical projection $G_S \twoheadrightarrow G$ and put $\Psi := \chi \circ \pi_K$, so that
$\Psi(\phi_\nu^{-1})=\chi^{-1}(\phi_{K,\nu})$ for all $\nu\not\in S$.
    Clearly $\Psi$ is a continuous character of $G_S$ and 
    \begin{eqnarray}
        \Psi \left( \Theta_S(X) \right)(q^{-s})    & = & \prod_{\nu\not\in S} \left(1 - \chi^{-1}(\phi_{K,\nu})
                                                    q^{-s d_\nu} \right)^{-1}  \nonumber\\
                                                & = & L_{K/F}(s, \chi^{-1}) \prod_{\nu \in S - S_K}\left(1 -
                                                    \chi^{-1}(\nu)q^{-s d_\nu} \right) \,. \nonumber
    \end{eqnarray}
(b) The profinite group $G_S / \Ker(\Psi)$ is topologically isomorphic to $\Psi (G_S)$ which is a
compact subgroup of $\C^{\times}$, so $\Ker(\Psi)$ has finite index. 
    Let $K$ be the fixed field of $\Ker(\Psi)$, so that $\Gal(K/F)\simeq G_S / \Ker(\Psi)$, and let $\chi$ be the character induced by $\Psi$
    on this quotient.  By definition $\Psi(\phi_\nu) = \chi(\phi_{K,\nu})$ and
    \begin{eqnarray}
        \Psi \left( \Theta_S(X) \right)    & = & \prod_{\nu\not\in S} \left(1 - \Psi(\phi_\nu^{-1})X^{d_\nu}
                                            \right)^{-1}
                                       = \prod_{\nu\not\in S} \left(1 - \chi^{-1}(\phi_{K,\nu})
                                            X^{d_\nu} \right)^{-1} \nonumber\\
                                        & = & \prod_{\nu\not\in S_K} \left(1 - \chi^{-1}(\nu)X^{d_\nu} \right)^{-1}
                                            \prod_{\nu\in S-S_K} \left(1 - \chi^{-1}(\nu)X^{d_\nu} \right) \,. \nonumber
    \end{eqnarray}
    Equation \eqref{rappresentazione complessa} follows immediately evaluating at $X=q^{-s}$.
\end{proof}

\subsection{Stickelberger series and the Tate algebra}\label{SecStickTate}
As an application of the previous theorem we prove that the Stickelberger series lies in the Tate algebra. \\
Let $\mathcal{R}$ be any topological ring. The Tate algebra $\mathcal{R}\langle X \rangle$ is the set
of formal power series whose coefficients converge to $0$, in
particular it contains the polynomial ring $\mathcal{R}[X]$.
Let $\mathcal{O}$ be the ring of integers of a finite extension of $\Qp$. We will be mainly interested in rings
of the form $\mathcal{R} = \mathcal{O} \llbracket \Gamma \rrbracket$, where $\Gamma$ will be the Galois group 
of an infinite extension of function fields. We recall that the topology on this ring is the weakest such that the projection 
$\mathcal{O} \llbracket \Gamma \rrbracket \twoheadrightarrow \mathcal{O} [\Gal(K/F)]$ is continuous for each finite subextension 
$K/F$: hence a sequence of elements $a_n$ of
$\mathcal{O} \llbracket \Gamma \rrbracket$ converges to $0$ if and only if the sequence of their projections is equal to
$0$ when $n$ is big enough, for each finite subextension $K/F$. \\
We consider $\Theta_S(X)$ as an element of $\mathcal{O}[G_S]\llbracket X \rrbracket$ via the (continuous) 
embedding $\Z \hookrightarrow \mathcal{O}$.

\begin{prop}\label{TateAlg} 
The series $\Theta_S(X)$ is an element of the Tate algebra $\mathcal{O} \llbracket G_S \rrbracket \langle X \rangle$.
\end{prop}

\begin{proof}
	Since $(1 - qX)^{-1} = \displaystyle{\sum_{n\geqslant 0}q^nX^n}$ is 
	a unit in $\mathcal{O} \llbracket G_S \rrbracket \langle X \rangle$, it is enough to show that $f(X) := (1 - qX)\Theta_S(X)$
	is in the Tate algebra. 
	
	\noindent Let $\Psi : G_S \rightarrow \C^\times$ be a continuous character. Following part (b) of Theorem \ref{ThmCollegamentoArtinL},
	let $K$ be the fixed field of $\Ker(\Psi)$ and $\chi$ the character induced by $\Psi$ on $\Gal(K/F)$ (note that, by definition, $\chi$ is trivial 
	if and only if $\Psi$ is trivial on $G_S$, which leads to $K=F$ and $S_F= \{\infty\}$).
	Hence we have
	\[
      \Psi \left( \Theta_S (X) \right) (q^{-s}) =
            L_{K/F}(s, \chi^{-1}) \prod_{\nu \in S - S_K}\left(1 - \chi^{-1}(\nu)q^{-s d_{\nu}} \right).
    \]
    The {\em full} Artin $L$-function is defined by the Euler product on all primes $\nu$
    \[
    \widetilde{L}_{K/F}(s, \chi^{-1}) = \prod_{\nu} \left(1 - \chi^{-1}(\nu) (\textbf{N}\nu)^{-s} \right)^{-1}
	\]
    and differs from $L_{K/F}(s, \chi^{-1})$ only for the factors associated to the primes of $S_K$. Thus 
    \begin{align*}
	    \Psi\left(f(X)\right)(q^{-s}) &= (1-q^{1-s})\Psi \left( \Theta_S (X) \right) (q^{-s}) \\
	    &= (1-q^{1-s})\widetilde{L}_{K/F}(s, \chi^{-1}) \prod_{\nu \in S}\left(1 - \chi^{-1}(\nu)q^{-s d_{\nu}} \right).
	\end{align*}
    A theorem of Weil \cite[V, Th\'eor\`em 2.5]{TateStark} shows that if $\chi \neq \chi_0$, then $\widetilde{L}_{K/F}(s, \chi^{-1})$
    is a polynomial in $q^{-s}$, and, for $\chi = \chi_0$, we have
    \[
    	\widetilde{L}_{K/F}(s, \chi_0)= \frac{P(q^{-s})}{(1-q^{-s})(1-q^{1-s})},
    \]
    where $P(X)$ is a polynomial of degree $2g-2$ ($g$ is the genus of $F$). Hence $\Psi\left(f(X)\right)$ is a polynomial
    for each continuous character $\Psi$ (here we use the fact that $S$ is not empty):
    for each finite subextension $K/F$
        \[   \pi_K\left( f(X) \right) \in \Z[\Gal(K/F)][X] \]
    (with $\pi_K$ the natural projection) and $f(X) \in \Z \llbracket G_S \rrbracket \langle X \rangle \subset 
    \mathcal{O} \llbracket G_S \rrbracket \langle X \rangle$.
\end{proof}

In Section \ref{SecStickClassGroups} we will need to evaluate the Stickelberger series $\Theta_S(X)$ at some element of 
$\mathcal{O}\llbracket G_S \rrbracket$, the previous proposition grants us that if we take $x$ in the unit disk 
$\{ x \in \mathcal{O}\llbracket G_S \rrbracket : |x| \leqslant 1 \}$, then
$\Theta_S(x)$ converges.

\section{Stickelberger series and Goss Zeta-function}\label{SecStickGossZeta}
We begin by giving a short description of the Goss Zeta-function and of its main properties (most notably its entireness): some
results are known, we give references where avaliable and refer the reader to \cite{CoscPhD} for more details.

\subsection{ The Goss Zeta-function}\label{Goss Zeta-function}
Let $\C_\infty$ be the completion of a fixed algebraic closure of $F_\infty$ and put $\Si:=\C_\infty^\times \times \Z_p$.
The analogue of the Dedekind Zeta-function for $F$ has been originally defined 
for some special values (the integers) by Carlitz in \cite{Carlitz} and later extended by Goss as a $\C_\infty$-valued function whose domain is $\Si$
 in \cite{Goss Zeta}. This work of Goss may be interpreted as a sort of analytic continuation of the function defined by Carlitz. 

\noindent We give the definition and a summary of the main properties of the {\em Goss Zeta-function}, details and proofs
can be found in \cite[Chapter 8]{Goss libro} or in \cite[Sections 1.4 and 1.5]{CoscPhD}.
  
\begin{defi}\label{DefSgn}
    A {\em sign function} on $F_\infty$ is any homomorphism ${\rm sgn} : F_{\infty}^{\times} \rightarrow
    \F_{\infty}^{\times}$ such that its restriction to $\F_{\infty}^{\times}$ is the identity.
    We extend ${\rm sgn}$ to all $F_\infty$ by defining ${\rm sgn}(0) = 0$.
\end{defi}

\noindent We fix a generator $\pi_\infty$ of the maximal ideal of $F_\infty$ and say that the
sign function sgn is {\em normalized} if $\text{sgn}(\pi_\infty) = 1$. Since $U_1(\infty)$ is a pro-$p$-group
and the image of sgn has order prime to $p$, every sign function is
trivial on $U_1(\infty)$. The isomorphism
\begin{equation}\label{decomposizione gruppo moltiplicativo}
    F_\infty^\times \simeq \pi_{\infty}^{\Z}\times \F_{\infty}^{\times} \times U_1(\infty)
\end{equation}
allows us to write any $a\in F_\infty^\times$ uniquely as
\begin{displaymath}
    a = \pi_{\infty}^{v_{\infty}(a)}\cdot \text{sgn}(a)\cdot \langle a \rangle_{\infty}. 
\end{displaymath}
We say an element $a\in F$ is {\em positive} if $\text{sgn}(a)=1$ and denote by $A_+$ the set of positive elements in $A$.
Let $\mathcal{I}$ be the set of nonzero fractional ideals of $F$ and denote by $\mathcal{P}_+$ the principal fractional ideals 
with a positive generator. The group $\mathcal{I}/\mathcal{P}_+$ is finite and we put $h^+(A):=|\mathcal{I}/\mathcal{P}_+|$:
it is easy to see that $h^+(A) = h^0(F) \cdot d_\infty \cdot (q^{d_\infty}-1) / (q-1)$.
We recall that $d_\infty=[\F_\infty :\F_q]$ and, for any 
$a\in F_\infty^\times$, define the {\em degree} of $a$ as $\deg(a) = - d_\infty v_{\infty}(a)$.
Note that if $I = (i)$ is principal, then the definition of $\deg(i)$ coincides with the degree of the ideal $I$,
i.e. $\deg(i)=\deg(I):=|A/I|$.


\noindent For every $u \in U_1(\infty)$ and $y \in \Z_p$ the series $\displaystyle{\sum_{n \geqslant 0} \binom{y}{n}(u-1)^{n}}$
converges in $U_1(\infty)$, so we put
\begin{displaymath}
    u^y = ((u-1)+1)^y := \sum_{n \geqslant 0} \binom{y}{n}(u-1)^{n} \,.
\end{displaymath}

\begin{defi}\label{DefIdealExp}
Given any ideal $I$, there exists a positive element $\alpha\in F^\times$ such that $I^{h^+(A)}=(\alpha)$. 
Put $\langle I\rangle_\infty:= \langle\alpha\rangle_\infty^{1/h^+(A)}$ as the unique $h^+(A)$-th root of $\langle\alpha\rangle_\infty$
which is still a $1$-unit. Then, for any $s=(x,y) \in \Si$, define 
\[ I^s:= x^{\deg(I)}\langle I\rangle_{\infty}^y \,.\]
\end{defi}

\noindent Fix a $d_\infty$-th root of $\pi_\infty$ and call it $\pi_*$: what follows will partly depend on this choice (but see
statement (f) below).
For every integer $j$ we put $s_j = (\pi_*^{-j}, j)$: the map $j \mapsto s_j$ gives us an embedding $\Z \hookrightarrow \Si$.
It is easy to see that $\langle I\rangle_\infty$ and $I^s$ are well defined and belong to a finite extension of $F$: 
we list a few fundamental properties of this exponential.

\begin{prop}\label{Exponential}
For every positive $a,b\in F_\infty^\times$, $s,t\in \Si$, $i,j\in \Z$ and every ideal $I\in \mathcal{I}$ one has
    \begin{enumerate}
    \item [(a)] $(ab)^s = a^s b^s$ and $a^{s+t} =a^s a^t$;
    \item [(b)] $a^{s_i} = a^i$ and $(a^{s_i})^{s_j} = a^{s_{ij}}$;
    \item [(c)] $I^{s_j}$ is algebraic over $F$ and, if $I=(\alpha) \in \mathcal{P}_+$, then $I^s = \alpha^s$;
    \item [(d)] let $F_\textbf{V}$ be the extension of $F$ obtained by adding every element of the form
            $I^{s_1}$, then $F_\textbf{V}/F$ is a finite extension with
            degree at most $h^+(A)$;
    \item [(e)] let $F_{\infty,\textbf{V}}$ be the extension of $F_\infty$ obtained by adding every
            element of the form $\langle I \rangle_{\infty}$, then
            $F_{\infty,\textbf{V}} / F_\infty$ is a finite $p$-extension with degree dividing $h^+(A)$;
    \item  [(f)] if $F$ contains all the $d_\infty$-th roots of unity, then $F_\textbf{V}$ does not depend on the
     choice of $\pi_\infty$ and $\pi_*$.
    \end{enumerate}
\end{prop}

\begin{proof} See \cite[Sections 8.1, 8.2]{Goss libro}, other details can be found in \cite[Section 1.4]{CoscPhD}.
\end{proof}

\begin{defi}
    The Goss Zeta-function is defined by the sum
    \begin{displaymath}
        \zeta_A (s) = \sum_{\mathfrak{a} \in \mathcal{I},\,\mathfrak{a} \geqslant 0} \mathfrak{a}^{-s}
        =\sum_{n \geqslant 0} \bigg(\sum_{\substack{\mathfrak{a} \in \mathcal{I},\, \mathfrak{a} \geqslant 0 \\ \deg(\mathfrak{a}) = n}}
                \langle \mathfrak{a} \rangle_{\infty}^{-y}\bigg)x^{-n} := \sum_{n \geqslant 0} a_n(y)x^{-n} 
    \end{displaymath}
    for $s=(x,y) \in \Si$.   
It converges for $|x|_\infty > 1$ and can also be rewritten as an Euler product 
\[ \zeta_A (s) = \prod_{\nu \neq \infty} \left(1 - \nu^{-s} \right)^{-1} =
	\prod_{\nu \neq \infty} \left(1 - \langle \nu\rangle_\infty^{-y} x^{-d_\nu)} \right)^{-1} \,. \]
\end{defi}

Estimates on the coefficients $a_n(y)$ allow to prove that the Goss Zeta-function is entire on $\Si$ (in the sense of \cite[Section 8.5]{Goss libro}), i.e.
the series provides an analytic continuation of $\zeta_A(s)$ to the whole $\Si$ (for a sketch of the proof see \cite[Sections 8.8 and 8.9]{Goss libro},
more details are in \cite[Section 1.5]{CoscPhD}).

\begin{teo}\label{ThmGossZetaExt}{\rm (Analytic extension of the Goss Zeta-function)}
	The serie $\zeta_A(s)$ is absolutely convergent for every $s= (x,y) \in \Si$ and is also uniformly convergent on the compact subsets
	of $\Si$.
\end{teo}

\subsection{Interpolation of the Goss Zeta-function}\label{SecInterpolation of the Goss Zeta-function}
Let $W_S$ be the subgroup of $G_S$ generated by all Artin symbols $\phi_\nu$ with $\nu \not\in S$ and
let $K$ be the fixed field of the topological closure of $W_S$. Since $\phi_\nu$ is a topological generator 
of the decomposition group of $\nu$ in $G_S$, the extension $K/F$ is totally split at every prime $\nu \not\in S$
and the Tchebotarev density theorem yields $K = F$, i.e. $G_S$ is the topological closure of
$W_S$.

\begin{lemma}\label{lemma frobenius diversi per primi diversi}
    Let $\lambda$ and $\mu$ be two distinct primes outside $S$, then $\phi_\lambda \neq \phi_\mu$.
\end{lemma}

\begin{proof}
    This is just class field theory: consider the following subgroup of the
    id\'eles $\If$ of $F$
    \begin{displaymath}
      H :=  F_\mu^\times \times \prod_{\nu \neq \mu ,\, \nu \not\in S} O_\nu^\times
                    \times \prod_{\nu\in S}\{ 1 \},
    \end{displaymath}
    and let $K$ be the class field of $F^\times H$ (as usual $F^\times$ is embedded diagonally
    in $\If$). Since $K/F$ is unramified outside $S$, totally split at $\mu$ and 
    inert in $\lambda$, the decomposition
    groups of $\mu$ and $\lambda$ in $G_S$ do not coincide and their Artin symbols are distinct.
\end{proof}

\noindent Let $f:=[F_{\infty,\textbf{V}}:F_\infty]$ (recall that $f$ is a power of $p$ by Proposition \ref{Exponential} part (e))
and let $N : F_{\infty,\textbf{V}}^\times \rightarrow F_\infty^\times$ be the norm map. For any $y \in \Z_p$ and
any $\nu \not\in S$ we define 
\begin{equation}\label{EqDefPsiy} 
\Psi_y (\phi_\nu) = N \left( \langle \nu \rangle_\infty^{-1} \right)^{y/f} \,.
\end{equation}
This is well defined because the norm sends $1$-units to $1$-units and it is possible to take the
$f$-th root without ambiguity.

\begin{lemma}
\label{lemma estension carattere psi}
    The map $\Psi_y$ extends to a group homomorphism $\Psi_y : G_S \rightarrow \C_\infty^\times$
    and induces a $\C_\infty^\times$-character on $\If$.
\end{lemma}

\begin{proof}
    Since $W_S$ is generated by the Artin symbols, we just put
    \[ \Psi_y (\sigma) = \Psi_y\bigg(\prod_{\nu\not\in S} \phi_\nu^{n_\nu}\bigg)=\prod_{\nu\not\in S} \Psi_y (\phi_\nu)^{n_\nu} .\] 
    Now $\Psi_y : W_S \rightarrow \C_\infty^\times$ is a continuous homomorphism (as composition of continuous maps) and it can be extended to the topological closure $G_S$.
    It is well defined (by Lemma \ref{lemma frobenius diversi per primi diversi}) and can be extended to id\'eles via the natural injective map
    \[ \varphi: \mathcal{I}\longrightarrow \If \,,\ {\rm given\ by}\quad
     \varphi(\mathfrak{a}) =\prod_\nu \pi_\nu^{v_\nu(\mathfrak{a})} ,\]
    and the reciprocity map $\text{rec}_S : \If \rightarrow G_S$. Since
    $\Ker (\text{rec}_S) = \displaystyle{ F^\times \cdot \prod_{\nu \not\in S}O_\nu^\times}$, 
    the map ${\rm rec}_S\circ\varphi$ is still injective. Obviously $\sigma=\displaystyle{\prod_{\nu\not\in S} \phi_\nu^{n_\nu}}\in W_S$
    is $({\rm rec}_S\circ \varphi)\bigg( \displaystyle{\prod_{\nu\not\in S} \nu^{n_\nu}}\bigg)=
    ({\rm rec}_S\circ \varphi)(\mathfrak{a})$ for some fractional ideal $\mathfrak{a}$ in
    $\mathcal{I}$, and $\Psi_y(\sigma)$ does not depend on the chosen expression for $\sigma$.
\end{proof}

\noindent We mention that the interpretation of $\Psi_y$ as a $\C_\infty^\times$-character on $S$-id\'eles
is the approach suggested in \cite[Theorem 3.8 and Remark 3.9]{ABBL}: we shall see a more
explicit relation between $\Psi_y$ and ${\rm rec}_S$ in the special case presented in the next section.

\noindent The extension $\Psi_y : \Z [G_S]\llbracket X \rrbracket \rightarrow \C_\infty \llbracket X \rrbracket$
gives an interpolation formula for the Goss Zeta-function (the case of the 
Carlitz cyclotomic extension is presented in \cite[Theorem 4.2]{BBC}).

\begin{teo}
\label{teorema stickelberger-goss}
    For every $s=(x,y) \in \Si$ we have
    \begin{displaymath}
        \Psi_y \left( \Theta_S (X) \right) (x) = \zeta_A (-s) \prod_{\nu \in S ,\, \nu \neq \infty}(1-\nu^s) \,.
    \end{displaymath}
\end{teo}

\begin{proof}
    Let $\nu$ be a prime outside $S$, assume $[\nu]$ has order $n$ in $\mathcal{I} / \mathcal{P}_+$ and let $\alpha$
    be a positive element such that $\nu^n = (\alpha)$. By definition $\langle \nu\rangle_\infty^n=\langle\alpha\rangle_\infty$.
    Now write $n=p^r n'$ with $(p, n')=1$ and let $u\in U_1(\infty)$ be the only
    $1$-unit verifying $u^{n'}=\langle \alpha \rangle_\infty$, so that $\langle \nu \rangle_\infty$ is
    a root of $f(X) = X^{p^r} - u \in F_\infty[X]$. Let $g(X)$ be the minimal polynomial of
    $\langle \nu \rangle_\infty$ over $F_\infty$. Since $f(X)$ is totally inseparable, it must be $f(X) = g(X)^{p^l}$ 
    and $g(X) = X^{p^k} - v$, where $l, k$ and $v$ satisfy $r= k+ l$ and $u = v^{p^l}$. Then $F_\infty(\langle \nu \rangle_\infty)/F_\infty$
    has degree $p^k$, while $F_{\infty,\textbf{V}} / F_\infty(\langle \nu \rangle_\infty)$ has degree $f / p^k$. Therefore
    \begin{displaymath}
        N\left(\langle \nu \rangle_\infty\right) = N_{F_\infty(\langle \nu \rangle_\infty)/F_\infty}(\langle \nu \rangle_\infty)^{f/p^k}
        = v^{f/p^k} = \langle \nu \rangle_\infty^f 
    \end{displaymath}
    (everything works for $p=2$ as well since in that case $v =- v$).\\
    Hence $\Psi_y (\phi_\nu) = \langle \nu \rangle_\infty^{-y}$ and 
    \[ \begin{array}{ll} \Psi_y \left( \Theta_S (X) \right)(x) & = \displaystyle{ \prod_{\nu\not\in S} \left(1 - \Psi_y(\phi_\nu^{-1})
                                                x^{d_\nu} \right)^{-1} = \prod_{\nu\not\in S} \left(1 - \langle \nu \rangle_\infty^{y}
                                                x^{d_\nu} \right)^{-1} }\nonumber \\
                                           & =  \displaystyle{\zeta_A(-s) \prod_{\nu\in S ,\, \nu\neq\infty} (1 - \nu^{s}) \,.}\hspace{4truecm}\nonumber \qedhere
    \end{array} \]
\end{proof}

\subsection{A special case} 
Here we provide a link between $\Psi_y$ and the Artin reciprocity map assuming that the 
class number of $F$ is $1$ and that $d_\infty=1$ (the hypotheses we shall use in Section \ref{SecStickClassGroups}): 
in particular $\Fi=\F_q$ and every element of
$F^\times$ can be written uniquely as product of a constant and a positive element. 
We choose an uniformizer at $\infty$ by simply taking a prime $\mathfrak{p} \neq \infty$ of degree $1$, 
letting $\pi_{\mathfrak{p}}$ be its unique positive generator and then putting $\pi_\infty := \pi_{\mathfrak{p}}^{-1}$. 
Note that this uniformizer is positive and is an element of $F$ \footnote{If $F$ is the rational function field
$\Fq(t)$, and $\infty=\frac{1}{t}$, we are simply taking as uniformizer the element $\pi_\infty = 1/(t-\alpha)$, 
where $\alpha$ is any element of $\Fq^\times$ and $\mathfrak{p} = (t - \alpha)$. }. Finally note that, in the
notations of Proposition \ref{Exponential}, $F_{\textbf{V}} = F$ and $F_{\infty,\textbf{V}} = F_\infty$:
therefore $\Psi_y(\phi_\nu) = \langle \pi_\nu \rangle_\infty^{-y}$ belongs to $U_1(\infty)$. 

\begin{teo} \label{ThmIsoH}
Let $F_{\infty, +}^\times$ be the kernel of the sign function. The natural inclusion 
    \begin{displaymath}
        F_{\infty, +}^\times \times \prod_{\nu \neq \infty} O_\nu^\times =: \mathcal{H} \hookrightarrow \If
    \end{displaymath}
    induces an isomorphism $\mathcal{H}\simeq \If / F^\times$.
    \end{teo}

\begin{proof}
    Let ${\bf i} = (i_\infty, i_{\nu_1}, i_{\nu_2}, \dots)$
    and ${\bf j} = (j_\infty, j_{\nu_1}, j_{\nu_2}, \dots)$ be two id\'eles in $\mathcal{H}$ with the same image in $\If/F^\times$, i.e. 
    there is an $x \in F^\times$ such that
    $\textbf{i} = x\textbf{j}$. For every $\nu\neq\infty$, $i_\nu = x j_\nu$ yields $v_\nu(x)=0$ because $i_\nu$ and $j_\nu$ are units.
    Moreover the product formula implies $v_\infty(x) = \displaystyle{-\sum_{\nu\neq\infty} d_\nu v_\nu(x) = 0}$ and so $x$ is a 
    constant. Finally since $i_\infty = x j_\infty$ and both $i_\infty$ and $j_\infty$ are positive, 
    we have $x=1$. \\
    To complete the proof we need surjectivity: let $\textbf{i} = (i_\infty, i_{\nu_1}, i_{\nu_2}, \dots)\in \If$ and put
    \begin{displaymath}
        x_{\bf i} = \text{sgn}(i_\infty)^{-1}\prod_{\nu \neq \infty}
                            \pi_\nu^{-v_\nu(i_\nu)} \in F^\times\,.
    \end{displaymath}
    It is easy to check that the id\'ele $x_{\bf i}\textbf{i}$ is in $\mathcal{H}$, and the proof is complete.
\end{proof}

Let $C_F:=\If /F^\times$ and consider the subgroup $O_S = \displaystyle{\prod_{\nu\not\in S}O_\nu^\times}/F^\times$: by
class field theory, the Artin map induces a continuous embedding
$\text{rec}_S : C_F / O_S \hookrightarrow G_S$,
which is not surjective because $G_S$ is profinite, while the quotient
$C_F / O_S$ is not. Let $\widehat{C_F / O_S}$ be the profinite completion of
$C_F / O_S$, then the map $\text{rec}_S$ extends to an isomorphism of topological groups
$ \widehat{\text{rec}}_S : \widehat{C_F / O_S} \stackrel{\sim}{\longrightarrow} G_S$ and, by Theorem \ref{ThmIsoH}, one has the isomorphism
\begin{displaymath}
    \widehat{C_F / O_S}\simeq \widehat{\langle \pi_\infty\rangle} \times U_1(\infty) \times
    \prod_{\substack{ \nu \in S \\ \nu \neq \infty}} O_\nu^\times
\end{displaymath}
where $\widehat{\langle \pi_\infty\rangle} \simeq \widehat{\Z}$. We denote by 
$\widehat{\pi}$ the canonical projection $\widehat{C_F / O_S} \twoheadrightarrow U_1(\infty)$.

\begin{teo}
    For every $y \in \Zp$ we have the following commutative diagram
	\[ \xymatrix{ G_S \ar[r]^{\Psi_y\ \ } \ar[d]_{\widehat{\rm rec}_S^{-1}} & U_1(\infty) \\
     \widehat{C_F/O_S} \ar@{->>}[r]_{\widehat{\pi}\ }  & U_1(\infty), \ar[u]_{y} }\]
    where $y$ denotes the raise-to-the-power $y$ map.
\end{teo}

\begin{proof}
    We prove the case $y=1$ first. Since $G_S$ is the topological closure of $W_S$, and all maps are continuous, 
    it is enough to show that $\Psi_1(\phi_\nu) = \widehat{\pi} \circ \widehat{rec}_S^{-1}(\phi_\nu)$ for every $\nu\not\in S$. \\
    The local Artin map $\text{rec}_\nu$ sends $\pi_\nu$
    to $\phi_\nu$, because the extension $F_S / F$ is unramified at $\nu$.
    Let $\textbf{i}_\nu$ be the id\'ele whose $\nu$-coordinate is equal to $\pi_\nu$ and whose $\mu$-coordinates
    for $\mu\neq \nu$ are all equal to $1$, and let $[\textbf{i}_\nu]	\in C_F$ be its equivalence class. Then 
\[ \phi_\nu = \text{rec}_S \left( [\textbf{i}_\nu] \right) = 
\text{rec}_S \left( [ \pi_\nu^{-1}\textbf{i}_\nu] \right) \] 
and, noting that $\pi_\nu^{-1}\textbf{i}_\nu\in\mathcal{H}$, we obtain 
    \begin{displaymath}
        \widehat{\pi} \circ \widehat{rec}_S^{-1}(\phi_\nu) = \langle \pi_\nu^{-1}\rangle_\infty = \Psi_1(\phi_\nu) \,.
    \end{displaymath}
    For a general $y \in \Zp$ note that, for any Artin symbol
    $\phi_\nu$, we have $\Psi_y (\phi_\nu) = \Psi_1 (\phi_\nu)^y$.
\end{proof}

\section{Stickelberger series and $\nu$-adic Zeta-functions}\label{SecSticknuadicL}
Fix a place $\nu$ different from $\infty$: in this section we consider a $\nu$-adic analogue of the Goss Zeta-function. 
Let $F_\nu$, $\C_\nu$, $\F_\nu$ and $\pi_\nu$ be the
$\nu$-adic versions of the objects defined for the place $\infty$. Fix an algebraic
closure $\overline{F}$ of $F$ and let $\iota : \overline{F} \hookrightarrow \C_\nu$ be
an $F$-embedding that is the identity on the compositum of the algebraic closure of $\F_q$ and
$F_{\textbf{V}}$. All the objects that we define here depend on $\iota$,
but we will omit this dependency to simplify notations. \\
The field $F_{\nu,\textbf{V}} := \iota(F_\textbf{V})F_\nu$ will play the role of $F_{\infty,\textbf{V}}$: as in Proposition \ref{Exponential}, 
one can show that $F_{\nu,\textbf{V}}$ is a $p$-extension of $F_\nu$ with degree dividing $h^+(A)$. 
The residue field of $F_{\nu,\textbf{V}}$ is still $\F_\nu$ and the cyclic group $\Z / |\F_\nu^\times|$ acts on its multiplicative subgroup. 
Let $\Sv = \C_\nu^\times \times \Z_p \times \Z / |\F_\nu^\times|$, which is a subgroup of the group of
$\C_\nu^\times$-valued characters on $F_\nu^\times$, and take $s_\nu = (x,y,j) \in \Sv$: we have to define
the exponential $I^{s_\nu} \in \C_\nu^\times$, for every fractional ideal $I \in \mathcal{I}$ coprime with $\nu$.  

\noindent For $s_1=(\pi_*^{-1},1)\in \Si$ the element $I^{s_1} \in F_\textbf{V}$ is a root of the polynomial
$X^{h^+(A)} - \alpha$, where $\alpha$ is the unique positive generator of $I^{h^+(A)}$, so the valuation at $\nu$ of $I^{s_1}$ is equal to zero. 
Hence $\iota (I^{s_1})$ is a unit in $F_{\nu,\textbf{V}}$ and can be written uniquely as a product
\begin{displaymath}
    \iota (I^{s_1}) = \omega (I) \langle I \rangle_\nu\,,
\end{displaymath}
for some $\omega (I) \in \F_\nu^\times$ and $\langle I \rangle_\nu$ a $1$-unit of $F_{\nu,\textbf{V}}$.
With the above notation, it is easy to check that the map $\omega : \mathcal{I}_\nu \rightarrow \F_\nu^\times$ sending $I$ to $\omega(I)$
is a group homomorphism on the group $\mathcal{I}_\nu$ of fractional ideals prime with $\nu$
(it is basically a {\em Teichm\"uller character}). 

\noindent Finally for $s_\nu = (x,y,j) \in \Sv$ and $I\in\mathcal{I}_\nu$ we define
\begin{displaymath}
    I^{s_\nu} = x^{\deg(I)} \omega(I)^j \langle I \rangle_\nu^y \,.
\end{displaymath}
We can embed $\Z$ in $\Sv$ via the map $j \in \Z \mapsto s_{\nu,j}=(1,j,j) \in \Sv$ and
one can show that this $\nu$-adic exponential satisfies properties analogous to the ones of
Proposition \ref{Exponential}.

\begin{prop}
For every $s_{\nu},t_{\nu} \in \Sv$, $I,J \in \mathcal{I}_\nu$ and $i,j\in \Z$, one has
    \begin{enumerate}
    \item[$\bullet$] $I^{s_{\nu}+t_{\nu}} =I^{s_{\nu}} I^{t_{\nu}}$ and $(IJ)^{s_{\nu}} = I^{s_{\nu}} J^{s_{\nu}}$.
    \item[$\bullet$] $(I^{s_{\nu,i}})^{s_{\nu,j}} = I^{s_{\nu, ij}}$ and $I^{s_{\nu,j}} = I^{s_j}$. 
    In particular $I^{s_{\nu,j}}$ is algebraic over $F$.
    \end{enumerate}
\end{prop}

\begin{defi}\label{DefGossnu}
The {\em $\nu$-adic Goss Zeta-function} is defined as
\begin{displaymath}
    \zeta_\nu(s_\nu) = \sum_{\mathfrak{a} \in \mathcal{I}_\nu,\,\mathfrak{a} \geqslant 0} \mathfrak{a}^{-s_\nu}
                     = \prod_{\mathfrak{p} \neq \nu,\infty} (1 - \mathfrak{p}^{-s_\nu})^{-1} \,.
\end{displaymath}
\end{defi}

\subsection{Interpolation of the $\nu$-adic Zeta-function}\label{SecInterpolvadicZeta}

We look for an analogue of Theorem \ref{teorema stickelberger-goss} for $\zeta_\nu(s_\nu)$:
this will be done only for primes $\nu \in S - \{\infty\}$, so for the rest of this section we take such a prime
$\nu$ assuming $\#S\geqslant 2$. 
Let $f_\nu$ be the degree of the extension $F_{\nu,\textbf{V}} / F_\nu$ and denote by 
$N_\nu : F_{\nu,\textbf{V}}^\times \rightarrow F_\nu^\times$ the norm map. Take $(y,j) \in \Z_p \times \Z/ |\F_\nu^\times|$ and,
 for $\mathfrak{p} \not\in S$, put 
\[ \Psi_{y,j}(\phi_\mathfrak{p}) = 
N_\nu (\langle \mathfrak{p}\rangle_\nu^{-1})^{y/f_\nu}\omega(\mathfrak{p})^{-j} \,,\] 
which is well defined by Lemma \ref{lemma frobenius diversi per primi diversi}. Adapting the proof of Lemma \ref{lemma estension carattere psi}
one obtains

\begin{lemma}
For every $(y,j) \in \Z_p \times \Z/ |\F_\nu^\times|$ the map $\Psi_{y,j}$ extends to a
continuous ring homomorphism $\Z[G_S] \llbracket X \rrbracket \rightarrow \C_\nu\llbracket X \rrbracket$.
\end{lemma}

Here is the relation between Stickelberger series and the $\nu$-adic Zeta-function. 

\begin{teo}\label{teorema interpolazione v-adica}
    For every $s_\nu=(x,y,j) \in \Sv$ (with $\nu\in S-\{\infty\}$) we have
    \begin{displaymath}
        \Psi_{y,j} \left( \Theta_S (X) \right) (x) = \zeta_\nu (-s_\nu) \prod_{\substack{\mathfrak{p} \in S \\
            \mathfrak{p} \neq \nu,\infty}}(1-\mathfrak{p}^{s_\nu}).
    \end{displaymath}
\end{teo}

\begin{proof}
    Le $\mathfrak{p}$ be a prime not in $S$, $n$ the exact order of $[\mathfrak{p}]$ in
    $\mathcal{I} / \mathcal{P}_+$ and $\alpha$ a positive element such that
    $\mathfrak{p}^n = (\alpha)$. We have $(\mathfrak{p}^n)^{s_1} = (\alpha)^{s_1} = \alpha$ and so
    \begin{displaymath}
        \omega(\mathfrak{p})^n \langle \mathfrak{p}\rangle_\nu^n =
        \iota (\mathfrak{p}^{s_1})^n = \iota (\alpha)^{s_1} = \alpha \,,
    \end{displaymath}
    which leads to $\langle \mathfrak{p}\rangle_\nu^n = \langle \alpha\rangle_\nu$. \\
    Proceeding like in Theorem \ref{teorema stickelberger-goss}, one finds that 
    $F_{\nu,\textbf{V}}/F_\nu(\langle\mathfrak{p}_\nu\rangle)$ has degree $f_\nu/p^k$ and
    \begin{displaymath}
        N_\nu\left(\langle \mathfrak{p} \rangle_\nu\right) = N_{K/F_\nu}(\langle \mathfrak{p}
        \rangle_\nu)^{f_\nu/p^k} = \langle \mathfrak{p} \rangle_\nu^{f_\nu} \,.
    \end{displaymath}
    From this we obtain that $\Psi_{y,j} (\phi_\mathfrak{p}) =
    \omega(\mathfrak{p})^{-j}\langle \mathfrak{p} \rangle_\nu^{-y}$ and
    \[ \Psi_{y,j} \left( \Theta_S (X) \right) = 
    \prod_{\mathfrak{p}\not\in S} \left(1 - \Psi_{y,j} (\phi_\mathfrak{p}^{-1})X^{d_\mathfrak{p}} \right)^{-1} = 
    \prod_{\mathfrak{p}\not\in S} \left(1 - \omega(\mathfrak{p})^j\langle \mathfrak{p}\rangle_\nu^{y} X^{d_\mathfrak{p}} \right)^{-1} \ .\]
    Hence
    \[ \Psi_{y,j} \left( \Theta_S (X) \right)(x) = \prod_{\mathfrak{p}\not\in S} \left(1 -  \omega(\mathfrak{p})^j \langle \mathfrak{p}
        \rangle_\nu^{y}   x^{d_\mathfrak{p}} \right)^{-1}=\zeta_\nu(-s_\nu) \prod_{\substack{\mathfrak{p}
                                                    \in S \\ \mathfrak{p}\neq\nu,\infty}}(1 -
                                                    \mathfrak{p}^{s_\nu}) \,.\qquad \qedhere \]
\end{proof}

\subsection{Interpolation via a $\nu$-adic $L$-function}\label{SecvadicLfuncInterpolation}
We investigate the values of $\nu$-adic Zeta-functions at integers, introducing
$\nu$-adic $L$-series to interpolate them. This is one of the main features of
Iwasawa theory, where $\mathfrak{p}$-adic $L$-functions (for $\mathfrak{p}$ a prime in a global field) usually represent the analytic counterpart of
characteristic ideals of Iwasawa modules.  

We shall use the following (see \cite[Lemma 8.8.1]{Goss libro} or \cite[Theorem 5.1.2]{Thakur}).

\begin{lemma}
\label{lemma Goss algebra lineare}
    Let $K$ be a function field with constant field $\F_q$, $v$ any normalized valuation on $K$ and $W\subset K$
    an $\F_q$-vector space with finite dimension. Assume that $v(w)>0$ for every $w \in W$.
   If $i$ is an integer with $0 \leqslant i < (q-1) \dim_{\F_q} W$, then for every $x \in K$ we have
            \[  \sum_{w \in W} (x+w)^i = 0 \,. \]
\end{lemma}

\noindent For each pair of non negative integers $j$ and $n$, we define
\[    S_n(j) = \sum_{\mathfrak{a} \geqslant 0,\,\deg \mathfrak{a} = n} \mathfrak{a}^{s_j}
\quad{\rm and}\quad 
    Z(X,j) = \sum_{n \geqslant 0} S_n(j)X^n \in F_{\textbf{V}}\llbracket X\rrbracket \,.\]

\begin{lemma}
    The $Z(X,j)$ are polynomials of degree less than or equal to 
    $d_\infty + 2g + \lfloor \frac{j}{q-1}\rfloor$.
\end{lemma}

\begin{proof}
    Fix a non negative integer $j$ and an $n> d_\infty + 2g  + \lfloor \frac{j}{q-1}\rfloor$.
    
    \noindent For any $\ell=1,\dots,h^+$, fix a representative $\mathfrak{a}_\ell\in \mathcal{I}$ for the class $C_\ell$ of 
    $\mathcal{I} / \mathcal{P}_+$. Then
    \begin{displaymath}
     S_n(j) =\sum_{j=1}^{h^+}  \sum_{\substack{\mathfrak{a} \geqslant 0,\, \mathfrak{a} \in C_\ell \\ \deg(\mathfrak{a}) = n }} \mathfrak{a}^{s_j} 
                = \sum_{j=1}^{h^+} \mathfrak{a}_\ell^{s_j}\cdot\sum_{\substack{\alpha \in F_+ ,\, 
                    \alpha \mathfrak{a}_\ell \geqslant 0 \\ \deg(\alpha) = n - \deg(\mathfrak{a}_\ell)}} \alpha^j :=
                    \sum_{j=1}^{h^+} S_n(C_\ell, j)\,.
    \end{displaymath}
    To prove $S_n(C_\ell,j)=0$, we put $n_\ell = n - \deg(\mathfrak{a}_\ell)$ and only consider the case $d_\infty\mid n_\ell$, otherwise 
    the claim is trivial. Let $\mathcal{X}:=\{ \alpha\in F_+\,:\, \deg(\alpha)=n_\ell\ {\rm and}\ \alpha \mathfrak{a}_\ell\geqslant 0\}$,
    denote by $D_1$ the divisor ${\rm Div}(\mathfrak{a}_\ell) +  (n_\ell / d_\infty - 1)\infty$ and by 
    $\mathcal{L}(D_1):=\{ \alpha\in F^\times\,:\, {\rm Div}(\alpha)+D_1\geqslant 0\}$ the associated Riemann-Roch space.
    The set $\mathcal{L}(D_1)$ acts freely by traslation on $\mathcal{X}$: write $\mathcal{X}$ as the union of $r$ orbits $\mathcal{X}_l$ and
    fix a representative $x_l \in \mathcal{X}_l$. Then 
    \begin{displaymath}
        S_n(C_\ell, j) = \mathfrak{a}_\ell^{s_j} \sum_{l=1}^r \sum_{u \in \mathcal{L}(D_1)} (u + x_l)^j =
                \mathfrak{a}_\ell^{s_j}\sum_{l=1}^r x_l^j \sum_{w \in x_l^{-1}\mathcal{L}(D_1)} (w + 1)^j \,.
    \end{displaymath}
    We have $v_\infty (x_l) = - n_\ell / d_\infty$ and $v_\infty (u) \geqslant 1 - n_\ell / d_\infty$.
    Therefore $v_\infty (w)$ is positive for every $w \in x_l^{-1}\mathcal{L}(D_1)$ and
    the vector space $x_l^{-1}\mathcal{L}(D_1)$ satisfies the hypotheses of
    Lemma \ref{lemma Goss algebra lineare}. Hence the inner sum is zero when
    $j < (q-1)\dim_{\F_q}\mathcal{L}(D_1) = (q-1)(n - d_\infty -g + 1)$ (the last equality comes from the Riemann-Roch theorem).
    \end{proof}
    
The polynomials $Z(X,j)$ are strictly related to the special values of the Goss Zeta-function since we have that
$Z(1,j) = \zeta_A(-s_j)$ for any $j\in \N$. We provide here a $\nu$-adic interpolation of this polynomials, which will be 
linked to the special values of the $\nu$-adic Goss Zeta-function.

\begin{defi}\label{DefLnu}
For any $y \in \Zp$ and $i \in \Z / | \F_\nu^*|$, the {\em $\nu$-adic $L$-series} is defined by
\begin{displaymath}
    L_\nu(X, y, \omega^i) = \sum_{n \geqslant 0} \bigg( \sum_{\substack{\mathfrak{a}\in\mathcal{I}_\nu,\,\mathfrak{a} \geqslant 0\\
            \deg(\mathfrak{a}) = n }}
            \omega (\mathfrak{a})^i \langle \mathfrak{a} \rangle_\nu^y \bigg) X^n \,.
\end{displaymath} 
\end{defi}

\begin{prop}
	For every $s_\nu = (x,y,i) \in \Sv$ we have
	\begin{equation} \label{EqLnuInterpolatesZnu}
		L_\nu (x, y , \omega^i) = \zeta_\nu (-s_\nu).
	\end{equation}
\end{prop}
\begin{proof}
	For a fractional ideal $\mathfrak{a}$ coprime with $\nu$, we have $\mathfrak{a}^{s_\nu} = 
	\omega (\mathfrak{a})^i \langle \mathfrak{a} \rangle_\nu^y x^{\text{deg} \, \mathfrak{a}}$, so
	\begin{align*} 
	L_\nu (x, y , \omega^i) & = \sum_{n \geqslant 0}\  \sum_{\substack{\mathfrak{a}\in\mathcal{I}_\nu,\,\mathfrak{a} \geqslant 0\\
								\deg(\mathfrak{a}) = n }} \omega (\mathfrak{a})^i \langle \mathfrak{a} \rangle_\nu^y 
								x^n  \\
		\ & = \sum_{n \geqslant 0} \  \sum_{\substack{\mathfrak{a}\in\mathcal{I}_\nu,\,\mathfrak{a} \geqslant 0\\
								\deg(\mathfrak{a}) = n }} \mathfrak{a}^{s_\nu}    = 
								\sum_{\mathfrak{a}\in\mathcal{I}_\nu,\,\mathfrak{a} \geqslant 0} \mathfrak{a}^{s_\nu} = \zeta_\nu (-s_\nu)\,. \qquad\qedhere
	\end{align*} 
\end{proof}

\noindent The previous proposition, Theorem \ref{teorema interpolazione v-adica} and Section \ref{SecStickTate}
show that the function $L_\nu(X,y,\omega^i)$ converges on the closed unit disc: since our application to class groups will 
only require specialization at $X=1$ we do not deal with its entireness here (which anyway can be proved just like the 
entireness of $\zeta_A$).

\noindent The relation with the polynomials $Z(X,j)$ for some particular values of $i$ and $j$ can be made more explicit: in particular the following 
theorem shows that, for some $i$ and $j$, the series $L_\nu(X, y, \omega^i)$ is actually a polynomial.
Computations will shift between $\infty$ and $\nu$ using formulas seen in Section 
\ref{SecInterpolation of the Goss Zeta-function}.

\begin{teo}\label{teorema interpolazione v-adica2}
	Assume that $\nu \in S-\{\infty\}$.
    \begin{enumerate}
    \item[(a)] Let $i$ and $j$ be two non negative integers, such that $i \equiv j \pmod{q^{d_\nu} -1}$,
        then
        \begin{displaymath}
            L_\nu(X, j, \omega^i) = Z(X,j)(1- \nu^{s_j}X^{d_\nu}) \,.
        \end{displaymath}
        In particular $L_\nu(X, j, \omega^i)$ in a polynomial.
    \item[(b)] For every $y \in \Zp$ we have 
        \begin{displaymath}
            L_\nu(X, y, \omega^i) \equiv Z(X,i) \pmod{\overline{\nu}} \,,
        \end{displaymath}
        where $\overline{\nu}$ denotes any prime of $F_{\textbf{V}}$ lying above $\nu$.
    \end{enumerate}
\end{teo}

\begin{proof}
   (a) Using the unique factorization of ideals in $A$ and the multiplicativity on $\mathcal{I}$ 
    of the maps $I\mapsto\langle I \rangle_\nu$ and $I\mapsto \omega(I)$
    one rewrites
    \begin{equation}\label{equazione con nu - parte prima}
        L_\nu(X, j, \omega^i) = \prod_{\mathfrak{p} \neq \nu, \infty} \left( 1 -
                \omega (\mathfrak{p})^i \langle \mathfrak{p} \rangle_\nu^j
                X^{d_\mathfrak{p}}\right)^{-1} \,.
    \end{equation}
    In the proof of Theorem \ref{teorema interpolazione v-adica} we have seen that this product is equal to
    \begin{displaymath}
        \Psi_{j,i}\left( \Theta_S(X) \right) \cdot \prod_{\substack{\mathfrak{p} \in S \\
            \mathfrak{p} \neq \infty , \nu}}\left( 1 - \omega (\mathfrak{p})^i
            \langle \mathfrak{p} \rangle_\nu^j X^{d_\mathfrak{p}}\right)^{-1} \,.
    \end{displaymath}
    Since $i \equiv j \pmod{q^{d_\nu} -1}$ and $\iota$ is the identity on $\mathfrak{p}^{s_1}\in F_{\textbf{V}}$, one has
    \begin{displaymath}
        \Psi_{j,i} (\phi_\mathfrak{p}^{-1}) = \omega (\mathfrak{p})^i \langle \mathfrak{p}
            \rangle_\nu^j = \omega (\mathfrak{p})^j \langle \mathfrak{p}
            \rangle_\nu^j = \iota \left( \mathfrak{p}^{s_1}\right)^j = \mathfrak{p}^{s_j} = \pi_*^{-jd_\mathfrak{p}}
            \langle \mathfrak{p} \rangle_\infty^j = \Psi_j (\phi_\mathfrak{p}^{-1})
            \pi_*^{-jd_\mathfrak{p}} \,.
    \end{displaymath}
    Therefore
    \begin{eqnarray}
        \Psi_{j,i} \left( \Theta_S (X) \right) & = & \prod_{\mathfrak{p} \not\in S} \left( 1 -
                                                    \Psi_{j,i} (\phi_\mathfrak{p}^{-1})
                                                    X^{d_\mathfrak{p}}\right)^{-1} \nonumber \\
                                            & = & \prod_{\mathfrak{p} \not\in S}  \left( 1 -
                                                    \Psi_j (\phi_\mathfrak{p}^{-1})\pi_*^{-jd_\mathfrak{p}}
                                                    X^{d_\mathfrak{p}}\right)^{-1}  = \Psi_j \left( \Theta_S (\pi_*^{-j}X) \right) \,,\nonumber
    \end{eqnarray}
   and 
    \begin{equation}\label{equazione intermedia con L-nu}
        L_\nu(X, j, \omega^i) = \Psi_j \left( \Theta_S (\pi_*^{-j}X) \right) \cdot
            \prod_{\substack{\mathfrak{p} \in S \\
            \mathfrak{p} \neq \infty , \nu}}\left( 1 - \mathfrak{p}^{s_j} X^{d_\mathfrak{p}}\right)^{-1} \,.
    \end{equation}
    The same arguments used to obtain \eqref{equazione con nu - parte prima}, yield
    $ Z(X,j) = \displaystyle{ \prod_{\mathfrak{p} \neq \infty} }\left( 1 - \langle \mathfrak{p} \rangle_\infty^j
                (\pi_*^{-j}X)^{d_\mathfrak{p}}\right)^{-1}$ which can be rewritten as
    \begin{displaymath}
         Z(X,j) = \Psi_j \left( \Theta_S (\pi_*^{-j}X) \right) \cdot
            \prod_{\substack{\mathfrak{p} \in S \\
            \mathfrak{p} \neq \infty}}\left( 1 - \mathfrak{p}^{s_j} X^{d_\mathfrak{p}}\right)^{-1} \,.
    \end{displaymath}
    Together with \eqref{equazione intermedia con L-nu}, this leads to
       $ L_\nu(X, j, \omega^i) = Z(X,j)(1- \nu^{s_j}X^{d_\nu})$.
       
       \noindent    (b)    For every ideal $\mathfrak{a}$ we have 
	$\langle \mathfrak{a} \rangle_\nu^y \equiv
    1 \equiv \langle \mathfrak{a} \rangle_\nu^i  \pmod{\overline{\nu}}$. 
	Hence
    \[\begin{array}{ll}
            L_\nu(X, y, \omega^i)   & = \displaystyle {\sum_{n \geqslant 0} \bigg( \sum_{\substack{\mathfrak{a}\in\mathcal{I}_\nu,\,\mathfrak{a} \geqslant 0\\
                                          \deg (\mathfrak{a}) = n }}
                                          \omega (\mathfrak{a})^i \langle \mathfrak{a}
                                          \rangle_\nu^y \bigg) X^n } \nonumber \\
                                    &\equiv \displaystyle{ \sum_{n \geqslant 0} \bigg( \sum_{\substack{\mathfrak{a}\in\mathcal{I}_\nu,\,\mathfrak{a} \geqslant 0\\
                                          \deg (\mathfrak{a}) = n }}
                                          \omega (\mathfrak{a})^i \langle \mathfrak{a}
                                          \rangle_\nu^i \bigg) X^n \pmod{\overline{\nu}} } = L_\nu(X, i, \omega^i) \pmod{\overline{\nu}} \nonumber \\
                                      &  = Z(X,i)(1- \nu^{s_i}X^{d_\nu}) \pmod{\overline{\nu}} \equiv  Z(X,i) \pmod{\overline{\nu}} \,. \nonumber
    \qquad\qedhere\end{array} \]
\end{proof}

\section{Stickelberger series and class groups}\label{SecStickClassGroups}

In this section we shall deal with the algebraic side of the theory: we study the $p$-part of class groups of degree zero divisors
of the finite subextension of a $\Z_p^\infty$-extension of $F$ generated by the torsion of a
Hayes module, the main result is the computation of the Fitting ideal of the inverse limit of such groups 
(Theorem \ref{TeoFittIdealClGrIwasawaModule}), which turns out to be generated by (a specialization of) our
Stickelberger series. This provides one instance of the Main Conjecture in our setting: the link
between the Fitting ideal (i.e. the Stickelberger series) and the various $L$-functions has been described in details in the
previous sections. Our strategy, which heavily relies on the computations of Greither and Popescu in \cite{Greither Popescu Galois} and
\cite{Greither Popescu Fitting}, puts emphasis on finite subextensions: this approach, closer to the classical one in characteristic zero, has been suggested and developed in \cite{ABBL} for $F=\F_q(t)$.

\begin{osse}\label{RemAltApproach} 
	An alternative path using limits of characteristic ideals has been described in \cite{BBL}, 
where the authors study $\Zp^\infty$-extensions of global function fields using $\Zp^d$-filtrations (and the Main
Conjecture for them proved originally by Crew in \cite{Crew} and more recently by Burns and Khue, Lai and Tan
in \cite{Bur}). This approach has been applied also to Iwasawa theory for elliptic curves and abelian varieties in \cite{BBL2}, 
building on structure theorems for Selmer groups (see e.g. \cite{Tan},
\cite{Tan2} and \cite{BV}) and the only avaliable Main Conjecture in this setting, i.e. the one for 
constant abelian varieties in \cite{LLTT}. Applications to a Main Conjecture for $\Z_p^\infty$-extensions and to Akashi series are 
provided in \cite[Theorem 3.10]{BBL2} and
\cite[Section 3]{BV2}. It would be interesting to try to apply the finite subextension
approach to this abelian varieties setting as well.
\end{osse} 

Some results of this section (basically Theorems \ref{Teorema fitting per caratteri di tipo 1 o 2} and
\ref{Teorema fitting Class Gr per caratteri di tipo 1 o 2})
already appeared in \cite[Section 3]{BBC} as a specialization of more general computations provided in \cite[Section 2]{BBC}.
We decided to include a short account here as well for completeness.  

\noindent Throughout this chapter we will assume that $d_\infty = 1$ and that $p$ does not divide $h^0(F)$.
Under these assumptions we have the following simplifications:
\begin{enumerate}
    \item[$\bullet$] the residue field $\Fi$ coincides with the field of constants $\Fq$;
    \item[$\bullet$] every principal ideal admits a positive generator;
    \item[$\bullet$] the class number of the ring of integers $A$ is equal to $h^0(F)$;
    \item[$\bullet$] the field $F_{\infty, \textbf{V}}$ coincides with the field $F_\infty$;
    \item[$\bullet$] for every $a \in F$: $\deg(a) = - v_\infty(a)$;
    \item[$\bullet$] $\pi_* = \pi_\infty$.
\end{enumerate}
As mentioned in the introduction the first assumption is needed to ensure that all extensions we deal
with are geometric and it is not really restrictive because we can reduce to this case by 
extending the constant field of $F$. The assumption on $h^0(F)$ is more crucial (but not too restrictive) because we are going 
to deal with the characters of a Galois group 
whose cardinality is divisible by $h^0(F)$ (see Section \ref{SecComplexChar}).

\subsection{Hayes extensions}
Let $H_A$ be the Hilbert class field of $A$, i.e. the maximal abelian unramified extension of $F$ which is 
totally split at $\infty$. Since the prime $\infty$ has degree $1$, 
the constant field of $H_A$ is $\Fq$. Class field theory implies that $Pic(A)\simeq \Gal(H_A / F)$ and the isomorphism 
is provided by the Artin reciprocity map: in particular the class of a fractional ideal $\mathfrak{a}$ is sent to its Frobenius 
in $\Gal(H_A/F)$ and, in case the support of $\mathfrak{a}$ is disjoint from $S$, this is simply the restriction of its 
Artin symbol $\phi_\mathfrak{a}\in G_S$.

\begin{defi} \label{DefHayesModule}
    We denote by $H_A\{\tau \}$ the ring of skew-polynomials in the variable $\tau$ with coefficients 
	in $H_A$, where $\tau f=f^q\tau$. A {\em Hayes module} (or {\em sign-normalized rank 1 Drinfeld module}) 
	is a homomorphism of $\Fq$-algebras 
	$\Phi : A \rightarrow H_A\{\tau \}$, such that:
    \begin{enumerate}
    \item[(a)] the image of $A$ is not contained in $H_A$;
    \item[(b)] for every $a \in A$ the coefficient of degree $0$ of $\Phi_a:=\Phi(a)$ is equal to $a$;
    \item[(c)] for every $a \in A$ the degree in $\tau$ of $\Phi_a$ is equal to $\deg(a)$ (i.e. $\Phi$ has rank 1);
    \item[(d)] for every $a \in A$, the leading coefficient of $\Phi_a$ is
        ${\rm sgn} (a)$ (i.e. $\Phi$ is sign-normalized).
    \end{enumerate}
\end{defi}

\noindent For details on the Hayes modules and on the properties mentioned below, the reader may refer to \cite[Chapter 7]{Goss libro},
\cite{Hayes Drinfeld} and \cite{Shu Kummer}.

For every $a \in A$ and $x \in \overline{F}$ put $a \cdot x := \Phi_a(x)$; this defines an $A$-module 
structure on $\overline{F}$. For any integral ideal $\mathfrak{a}$ of $A$, consider the left ideal of
$H_A\{\tau \}$ generated by all the elements $\Phi_a$ with $a \in \mathfrak{a}$; since $H_A\{\tau \}$
is right-euclidean, we have that every left ideal is principal and we denote by $\Phi_\mathfrak{a}$
the unique monic generator of the ideal $(\Phi_a\,:\,a\in\mathfrak{a})$. 

\begin{defi}\label{DefHayesTor}
The {\em $\mathfrak{a}$-torsion} of $\overline{F}$ is the set
\[ \Phi[\mathfrak{a}]:=\{x\in \overline{F}\,:\,\Phi_\mathfrak{a}(x)=0 \} .\]
It is finite for any $\mathfrak{a}\neq 0$ and it is an $A/\mathfrak{a}$-module isomorphic to $A / \mathfrak{a}$.
\end{defi}

\noindent We put $F(\mathfrak{a}):=H_A(\Phi[\mathfrak{a}])$; the following theorem explains how 
the $\mathfrak{a}$-torsion of a Hayes module can be used to define extensions of $F$ analogous
to the cyclotomic extension of $\Q$. 

\begin{teo}\label{teorema riassuntivo sui gr galois estensioni ciclotomiche}
    The field $F(\mathfrak{a})$ is a geometric, abelian Galois extension of $F$ which verifies
	\begin{enumerate}
    \item[(a)] the only ramified primes in $F(\mathfrak{a}) / H_A$ are the primes of $H_A$ dividing
        $\mathfrak{a}$ and $\infty$;
    \item[(b)] $Gal \left(F(\mathfrak{a}) / H_A \right) \simeq (A / \mathfrak{a} )^\times$ via an isomorphism
        sending $a \in A$ to $\sigma_a$, where $\sigma_a(\lambda) = \Phi_a (\lambda)$ for every $\lambda \in \Phi[\mathfrak{a}]$;
    \item[(c)] the isomorphism in {\rm (b)} sends the decomposition and inertia groups of $\infty$  
        to $\Fq^\times$;
    \item[(d)] if $\mathfrak{p}^n$ is the exact power of $\mathfrak{p}$ dividing $\mathfrak{a}$,
        then the isomorphism in {\rm (b)} sends the inertia group of $\mathfrak{p}$ to $(A / \mathfrak{p}^n )^\times$.
    \end{enumerate}
\end{teo}

\noindent We fix a prime $\mathfrak{p}$ and put $S = \{ \mathfrak{p}, \infty \}$. For any $n\geqslant 0$ let $F_n := F(\mathfrak{p}^{n+1})$ 
and $G_n = {\rm Gal} (F_n / F)$. From part (b) of Theorem \ref{teorema riassuntivo sui gr galois estensioni ciclotomiche}
we have that $F_n / F_0$ is a finite $p$-extension and that $F_0 / F$ has degree $h^0(F)(q^{d_\mathfrak{p}} -1)$. 
Since we assumed $h^0(F)$ coprime with $p$, we obtain a decomposition $G_n\simeq G_0 \times \Gamma_n$, where 
$\Gamma_n = \Gal(F_n / F_0)$ is a $p$-group and $G_0$ has order prime with $p$. 
The fields $F_n$ form an Iwasawa tower: if we put $\mathcal{F}=\cup F_n$, then
\begin{displaymath}
    G_\infty := \Gal(\mathcal{F} / F)= \lim_{\leftarrow} G_n\simeq G_0 \times \lim_{\leftarrow} \Gamma_n
    =:G_0\times \Gamma_\infty \,,
    \end{displaymath}
with $\Gamma_\infty \simeq \Zp^\infty$. Note that the only primes which ramify in $\mathcal{F} / F$
are $\mathfrak{p}$ and $\infty$, so $\mathcal{F} \subseteq F_S$. The following
diagram gives a recap of the fields and Galois groups introduced above.
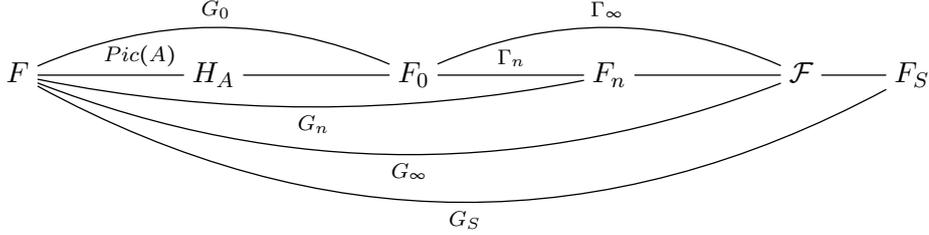
\[ \xymatrix{ F \ar@/^1.5pc/@{-}[rrrr]^{G_0} \ar@{-}[rr]^{\ \ \ \ \ Pic(A)} \ar@/_1pc/@{-}[rrrrrr]_{G_n} \ar@/_2.5pc/@{-}[rrrrrrrr]_{G_\infty}
\ar@/_4pc/@{-}[rrrrrrrrr]_{G_S} & & H_A \ar@{-}[rr] & & F_0 \ar@{-}[rr]^{\Gamma_n} \ar@/^1.5pc/@{-}[rrrr]^{\Gamma_\infty} & &
F_n \ar@{-}[rr] & & \mathcal{F} \ar@{-}[r] & F_S } \]

\noindent Regarding the behaviour of primes in these extensions, we recall that
\begin{itemize}
\item any prime different from $\mathfrak{p}$ and $\infty$ is unramified everywhere;
\item $\mathfrak{p}$ is unramified in $H_A / F$ and totally ramified in
$\mathcal{F} / H_A$;
\item $\infty$ is totally split in $H_A / F$, it ramifies (not totally in general) in $F_0 / H_A$ 
with inertia group isomorphic to $\Fq^\times$, and it is again totally split in $\mathcal{F}/F_0$.
\end{itemize}

\subsection{Fitting ideals and the Greither-Popescu theorem}\label{SecFitt1}
Fix an algebraic closure $\Falg$ of $\Fq$ and let $\gamma$ be the {\em arithmetic Frobenius}, i.e. a topological generator of 
$G_\F:=\Gal(\Falg / \Fq)$. For every field $L$ we put $L^{ar}:=\Falg L$; when $L=F_n$ (resp. $\mathcal{F}$),
we have that $F_n^{ar}$ (resp. $\mathcal{F}^{ar}$) is Galois over $F$, with Galois group isomorphic to $G_n \times G_\F$ 
(resp. $G_\infty \times G_\F$) since $F_n$ (resp. $\mathcal{F}$) is a geometric extension of  $F$. 

We recall here the definition of Fitting ideal of a finitely generated module. For an in-depth discussion the reader may refer to 
\cite[Chapter 3]{Northcott} or to the appendix of \cite{Mazur Wiles} for the main properties.
Let $R$ be any commutative and unitary ring, let $M$ be a finitely generated $R$-module and fix a set of generators 
$\{m_1,\dots, m_r\}$ for $M$. 
A \emph{relation vector} between the generators $m_i$ is an element $\underline{a} = (a_1,\dots,a_r) \in R^r$
such that $\sum a_im_i =0$. A \emph{matrix of relations} is any $q \times r$ matrix, with $q \geqslant r$, whose rows are 
relation vectors.

\begin{defi}
	The Fitting ideal of $M$, denoted ${\rm Fitt}_R(M)$, is the ideal generated by the determinants
	of all the $r\times r$ minors of all the matrices of relations of $M$. 
\end{defi}

\noindent It is well known that ${\rm Fitt}_R(M)$ does not depend on the chosen set of generators, moreover it is enough to 
consider square relation matrices of dimension $r$. We recall that ${\rm Fitt}_R(M)\subseteq {\rm Ann}_R(M)$,
i.e. we can have a nontrivial ideal only for torsion modules.

\subsubsection{The modules $H_n(\nu)$}
For every prime $\nu$ of $F$ there exists only a finite number of primes of $F_n^{ar}$ lying
above $\nu$, in particular
	\begin{itemize}
	\item there are $d_\mathfrak{p}  h^0(F)$ primes dividing $\mathfrak{p}$;
	\item there are $h^0(F)  q^{nd_\mathfrak{p}}  (q^{d_\mathfrak{p}}-1)/(q-1)$ primes dividing $\infty$.
	\end{itemize}

\noindent Let $H_n(\nu)$ be the free $\Zp$-module generated by the primes of $F_n^{ar}$
lying above $\nu$. Let $I_n(\nu) \subseteq G_n$ be the inertia group of $\nu$, since for $\nu=\infty$ such group does not depend on
$n$, we shall simply write $I_n(\infty)=:I_\infty$. We have that $H_n(\nu)$ is also a free $\Zp[G_n / I_n(\nu)]$-module of rank $d_\nu$, 
and there is a natural action of $G_\F$ on $H_n(\nu)$. We are interested in studying the structure of $H_n(\nu)$
as a $\Zp[G_n]\llbracket G_\F \rrbracket$-module. 

\noindent For $\nu \in \{ \mathfrak{p}, \infty \}$ we denote by ${\rm Fr}_\nu$ any lift to $G_n$ of the Frobenius map in $G_n / I_n(\nu)$. 
Since the decomposition and inertia groups of $\infty$ coincide,  we can choose ${\rm Fr}_\infty = 1$. The same choice is valid 
for $\mathfrak{p}$ if and only if $\mathfrak{p}$ is totally split in $H_A$. 

\begin{defi}\label{DefEulerFactor}
The {\em Euler factor} at $\nu$ is
\begin{displaymath}
    e_\nu(X) := 1 - {\rm Fr}_\nu^{-1} X^{d_\nu} \in \Zp[G_n]\llbracket X \rrbracket \,.
\end{displaymath}
Since we will also need to specialize the variable $X$ at $\gamma^{-1}$, we put
\begin{displaymath}
e_\nu := e_\nu(\gamma^{-1}) = 1 - {\rm Fr}_\nu^{-1} \gamma^{-d_\nu} \in \Zp[G_n]\llbracket G_\F \rrbracket \,.
\end{displaymath}
\end{defi}

\noindent The statements of \cite[Lemmas 2.1 and 2.2]{Greither Popescu Fitting}, adapted to our setting, 
translate into 

\begin{lemma}\label{lemma Greither-Popescu}
    For $\nu \in \{\mathfrak{p}, \infty\}$, let $Aug_{\nu,n}:=(\tau - 1\,:\,\tau \in I_n(\nu))$ be the {\em augumentation ideal}
    of $I_n(\nu)$ in $\Zp[G_n]\llbracket G_\F \rrbracket$. Then
    \begin{enumerate}
    \item[(a)] ${\rm Fitt}_{\Zp[G_n]\llbracket G_\F \rrbracket} \left(H_n(\nu)\right) = (e_\nu)$ for any $\nu \neq \mathfrak{p}, \infty$;
    \item[(b)] ${\rm Fitt}_{\Zp[G_n]\llbracket G_\F \rrbracket} \left(H_n(\infty)\right) =
            (e_\infty, Aug_{\infty, n})$;
    \item[(c)] ${\rm Fitt}_{\Zp[G_n]\llbracket G_\F \rrbracket} \left(H_n(\mathfrak{p})\right) =
            (e_\mathfrak{p}, Aug_{\mathfrak{p}, n})$.
    \end{enumerate}
\end{lemma}

\noindent Since the $H_n(\nu)$ are free $\Z_p[G_n/I_n(\nu)]$-modules, we have isomorphisms
\begin{enumerate}
\item[$\bullet$] if $\nu \neq \mathfrak{p}, \infty$
        \begin{displaymath}
            H_n(\nu) \simeq \Zp[G_n]\llbracket G_\F \rrbracket / (e_\nu) \,;
        \end{displaymath}
\item[$\bullet$] if $\nu = \infty$
        \begin{displaymath}
            H_n(\infty) \simeq \Zp[G_0 / \Fq^{\times} \times \Gamma_n]\llbracket G_\F \rrbracket / (e_\infty) \,;
        \end{displaymath}
\item[$\bullet$] if $\nu = \mathfrak{p}$
        \begin{displaymath}
            H_n(\mathfrak{p}) \simeq \Zp[Pic(A)]\llbracket G_\F \rrbracket / (e_\mathfrak{p}) \,.
        \end{displaymath}
\end{enumerate}

\subsubsection{Complex characters}\label{SecComplexChar}
Let $\chi \in \widehat{G_0}:=\text{Hom}(G_0, \C^\times)$ be a complex character for $G_0$, to include its values we shall extend scalars to
the Witt ring $W = \Zp[\zeta]$, where $\zeta$ denotes any primitive root of unity of order $|G_0|$ (recall that we are assuming $(|G_0|,p)=1$).
We let
\begin{displaymath}
    e_\chi := \frac{1}{|G_0|}\sum_{g \in G_0} \chi (g^{-1})g \in W[G_0]
\end{displaymath}
be the idempotent associated to $\chi$. For any $W[G_0]$-module $M$, we denote its $\chi$-part by $M(\chi):=e_\chi M$. 
For any $\Zp[G_0]$-module $M$ we shall abuse notations a bit and write $M(\chi)$ to denote $(W\otimes_{\Z_p} M)(\chi)$. 
This will not have any effect on the numerous exact sequences we are going to consider, because $W$ is a flat $\Zp$-module. 
Finally recall that if $G_0$ acts trivially on $M$, then 
\begin{displaymath}
M(\chi) = \left\{ \begin{array}{ll}
                W \otimes_{\Zp} M       &       \ \mbox{if } \chi=\chi_0 \,, \\
                \ & \\
                0                       &       \ \mbox{if } \chi \neq \chi_0 \,
    \end{array} \right.
\end{displaymath}
(as usual $\chi_0$ denotes the trivial character).

\begin{defi}\label{DefCharType}
    Let $\chi$ be a character of $G_0$. We will distinguish $3$ types of characters:
    \begin{enumerate}
        \item[$\bullet$] $\chi$ is said to be of \emph{type 1} if $\chi \left( I_\infty \right) \neq 1$;
        \item[$\bullet$] $\chi$ is said to be of \emph{type 2} if $\chi \left( I_\infty \right) = 1$ and
                    $\chi \left( \Gal(F_0 / H_A) \right) \neq 1$;
        \item[$\bullet$] $\chi$ is said to be of \emph{type 3} if $\chi \left( \Gal(F_0 / H_A) \right) = 1$.
    \end{enumerate}
\end{defi}

\noindent Recall that $d_\infty = 1$ and ${\rm Fr}_\infty =1$, so that $e_\infty = 1 - \gamma^{-1}$. Taking $\chi$-parts in 
Lemma \ref{lemma Greither-Popescu} we get
\begin{equation}\label{computo H-infinito}
H_n(\infty)(\chi) \simeq \left\{ \begin{array}{ll}
                            0               &       \ \mbox{if } \chi \mbox{ is of type } 1\,, \\
						\ & \\
						W[\Gamma_n]     &       \ \mbox{otherwise\,.}
                         \end{array} \right.
\end{equation}
For the prime $\mathfrak{p}$ we have the exact sequence (again from Lemma \ref{lemma Greither-Popescu})
\begin{equation}
	\xymatrix{e_\chi(1-{\rm Fr}_\mathfrak{p}^{-1}\gamma^{-d_\mathfrak{p}})W[Pic(A)]\llbracket G_\F \rrbracket \ar@{^(->}[r] & 
	  e_\chi W[Pic(A)]\llbracket G_\F \rrbracket \ar@{->>}[r] & H_n (\mathfrak{p})(\chi)\,,} 
	  \end{equation}
which leads to
\begin{equation}\label{computo H-p}
H_n(\mathfrak{p})(\chi) \simeq \left\{ \begin{array}{ll}
                            0                                   &       \ \mbox{if } \chi \mbox{ is of type }
                                                                            1 \mbox{ or } 2\,, \\
                            \ & \\
							W\llbracket G_\F \rrbracket / \left(1 - \chi({\rm Fr}_\mathfrak{p}^{-1})
                            \gamma^{-d_\mathfrak{p}}\right)           &       \ \mbox{otherwise\,.}
                         \end{array} \right.
\end{equation}

\subsubsection{The theorem of Greither and Popescu}
Let $X_n$ be the projective curve defined over $\mathbb{F}_q$ and associated with $F_n$; let 
\[ Jac(X_n)(\overline{\F}_q) \simeq {\rm Div}^0(F_n^{ar})/\{{\rm Div}(x)\,:\,x\in (F_n^{ar})^\times \} \]
be the set of $\overline{\F}_q$-rational points of its Jacobian.
Following \cite[Section 2]{Greither Popescu Galois},
we define the Deligne's Picard $1$-motive $\mathcal{M}_n:=\mathcal{M}_S(F_n^{ar})$ as the group
morphism ${\rm Div}^0(F_n^{ar})\rightarrow  Jac(X_n)(\overline{\F}_q)$
which induces the isomorphism above. In terms of \cite[Definition 2.3]{Greither Popescu Galois}
it is the morphism associated to $(X_n,\overline{\mathbb{F}}_q, S(F_n^{ar}),\emptyset)$, where $S(F_n^{ar})$ is
the set of primes of $F_n^{ar}$ lying above primes in $S$ and the choice of $\emptyset$ is justified by
\cite[Remark 2.7]{Greither Popescu Galois}. The $m$-torsion $\mathcal{M}_n[m]$ of $\mathcal{M}_n$ (see 
\cite[Definition 2.5]{Greither Popescu Galois}) fits into the exact sequence
\[ \xymatrix{ 0\ar[r] &  Jac(X_n)(\overline{\F}_q)[m] \ar[r] & \mathcal{M}_n[m] \ar[r] & 
{\rm Div}^0(S(F_n^{ar}))\otimes \Z/m\Z \ar[r] & 0 } \]
(where ${\rm Div}^0(S(F_n^{ar}))$ are the divisors supported on $S(F_n^{ar})$) and behaves well with respect to norm maps
so that we can define its $p$-adic Tate module as
$\displaystyle{ T_p(\mathcal{M}_n) = \lim_\leftarrow \mathcal{M}_n [p^m] }$.

\noindent We denote by $T_p(F_n) := T_p \left(Jac(X_n)(\overline{\F}_q) \right)$ the $p$-adic Tate module of 
the Jacobian of $X_n$. Our next task is to study the structure of $T_p(F_n)$ as a Galois module over
$\Zp[\Gamma_n]\llbracket G_\F \rrbracket$ and a crucial role will be played by Theorem \ref{Teo GP su ideale fitting del Deligne-Picard}
below.

\noindent Let $\Theta_n(X)$ (resp. $\Theta_\infty (X)$) be the projection of the Stickelberger series
$\Theta_S(X)$ to the ring $\Z [G_n] \llbracket X \rrbracket$ (resp. $\Z \llbracket G_\infty \rrbracket
\llbracket X \rrbracket$), which is easily seen to be the Stickelberger series associated to the
extension $F_n / F$ (resp. $\mathcal{F} / F$): actually one obtains an equivalent definition for $\Theta_S(X)$ by
taking the inverse limit (with respect to projections) of the Stickelberger series of the subextensions of $F_S$.

In \cite[Theorem 4.3]{Greither Popescu Galois} the authors prove the following

\begin{teo}\label{Teo GP su ideale fitting del Deligne-Picard}
One has ${\rm Fitt}_{\Zp[G_n]\llbracket G_\F \rrbracket} \left( T_p(\mathcal{M}_n)\right) =
        \left(\Theta_n(\gamma^{-1})\right)$.
\end{teo}

\noindent Note that evaluating the Stickelberger series $\Theta_n(X)$ at $X=\gamma^{-1}$ makes sense because of
Proposition \ref{TateAlg}.

\subsection{Fitting ideals for Tate modules: finite level}
We define $D_n$ to be the kernel of the $\deg$ map in the following exact sequence
\begin{displaymath}
 \xymatrix{ 0\ar[r] &  D_n \ar[r] & H_n(\infty) \oplus H_n(\mathfrak{p}) \ar[r]^{\qquad\ \ \deg} & \Zp \ar[r] & 0\,.}
\end{displaymath}
Since $G_0$ acts trivially on $\Zp$ we have
\begin{equation}\label{computo L_n}
D_n(\chi) \simeq \left\{ \begin{array}{ll}
                            0                                   &       \ \mbox{if } \chi \mbox{ is of type }
                                                                            1\,, \\
                            \ & \\                                                
                            W[\Gamma_n]                         &       \ \mbox{if } \chi \mbox{ is of type }
                                                                            2\,, \\
                            \ & \\                                                
                            W[\Gamma_n] \oplus W\llbracket G_\F \rrbracket / \left(1 -
                            \chi({\rm Fr}_\mathfrak{p}^{-1})\gamma^{-d_\mathfrak{p}}\right)           &       \ \mbox{if } \chi \mbox{ is of type }
                                                                            3 \\
                                                                            \ & \mbox{ and } \chi \neq \chi_0\,.
                         \end{array} \right.
\end{equation}
In \cite[after Definition 2.6]{Greither Popescu Galois} the authors provide the following exact sequence
\begin{equation}\label{successione esatta Greither Popescu}
\xymatrix{0\ar[r] & T_p(F_n) \ar[r] & T_p(\mathcal{M}_n) \ar[r] & D_n \ar[r] & 0 \,.}
\end{equation}
For every character $\chi$ we denote by $\Theta_n (X, \chi)$ the only element of
$\Z[\Gamma_n] \llbracket X \rrbracket$ that satisfies $\Theta_n(X,\chi)e_\chi =e_\chi \Theta_n(X)$. Then we have the following
(see also \cite[Theorem 3.2]{BBC})

\begin{teo}\label{Teorema fitting per caratteri di tipo 1 o 2}
    Let $\chi \in \widehat{G_0}$ be a character of type 1 or 2. Then
    \[
    	 {\rm Fitt}_{W[\Gamma_n]\llbracket G_\F \rrbracket} \left( T_p(F_n)(\chi)\right) =
    	 \left( \Theta_n^{\sharp}(\gamma^{-1},\chi) \right)\,,
    \]
    where we put    
    \begin{displaymath}
        \Theta_n^{\sharp}(\gamma^{-1},\chi) =
            \left\{ \begin{array}{ll}
                \Theta_n(\gamma^{-1},\chi)                        &   \ \mbox{if } \chi \mbox{ is of
                                                                                    type } 1\,, \\
                 \ & \\                                                                   
                \displaystyle{\frac{\Theta_n(\gamma^{-1},\chi)}{1-\gamma^{-1}} }  &   \ \mbox{if } \chi \mbox{ is of
                                                                                    type } 2\,.
            \end{array} \right.
    \end{displaymath}
\end{teo}

\begin{proof} Take $\chi$-parts in \eqref{successione esatta Greither Popescu} and use \eqref{computo L_n} to get:
\begin{itemize}
\item if $\chi$ is of type $1$, $D_n(\chi) = 0$ and $T_p(F_n)(\chi) \simeq T_p(\mathcal{M}_n)(\chi)$;
\item if $\chi$ is of type $2$, $D_n(\chi) = W[\Gamma_n]
    \simeq W[\Gamma_n]\llbracket G_\F \rrbracket /(1-\gamma^{-1})$ is a cyclic
    $W[\Gamma_n]\llbracket G_\F \rrbracket$-module: hence we can apply \cite[Lemma 3]{Cornacchia Greither}
    and obtain
    \begin{displaymath}
        (1-\gamma^{-1}) {\rm Fitt}_{W[\Gamma_n]\llbracket G_\F \rrbracket} \left( T_p(F_n)(\chi)\right)
           = {\rm Fitt}_{W[\Gamma_n]\llbracket G_\F \rrbracket} \left( T_p(\mathcal{M}_n)(\chi)\right) = \left(\Theta_n(\gamma^{-1},\chi)\right) . \qedhere
    \end{displaymath}
\end{itemize}
\end{proof}

When $\chi$ is a character of type $3$ things get more involved, since $D_n(\chi)$ is not cyclic
and we can only compute the Fitting ideal of some dual of the Tate module: details can be found in 
\cite[Section 2.5]{CoscPhD} or in \cite[Section 3]{BBC}.

\subsection{Fitting ideals for Tate modules: infinite level}\label{SubSecTateModInfLev}
To prove an analog of Theorem \ref{Teorema fitting per caratteri di tipo 1 o 2} for the {\em infinite level} $\mathcal{F}^{ar}$
we need to study the relation between $T_p(F_n)$ and $T_p(F_m)$ for $n\geqslant m$.
Let $\Gamma^n_m := \Gal(F_n / F_m)$ and $I_{\Gamma_m^n}$ the associated augmentation ideal in 
$W[\Gamma_n]$. We recall that $F_n/F_m$ is totally
ramified at all primes lying above $\mathfrak{p}$ and unramified everywhere else,
moreover the number of primes in $F_m$ above $\mathfrak{p}$ is the same for all $m$ and coincides
with the number of primes of $H_A$ lying above $\mathfrak{p}$.

\noindent We denote by $\overline{C}_n$ the $p$-part of the class group of degree zero divisors of
$F_n^{ar}$ and we recall that $T_p(F_n) = \text{Hom}(\Qp / \Zp , \overline{C}_n)$. Thus
there are natural maps induced by norms and inclusions respectively $N^n_m : T_p(F_n) \rightarrow T_p(F_m)$ and 
$i^m_n : T_p(F_m) \rightarrow T_p(F_n)$. We define
\begin{displaymath}
    T_p(\mathcal{F})(\chi) = \lim_\leftarrow T_p(F_n)(\chi) \,,
\end{displaymath}
where the limit is taken with respect to the norm maps. The limit
$T_p(\mathcal{F})(\chi)$ is a module over the profinite (non-noetherian) algebra $\Lambda_\F:=
W\llbracket \Gamma_\infty \rrbracket \llbracket G_\F \rrbracket$; our next goal is to compute its Fitting ideal. 

Norms and inclusions are defined on the modules $T_p(\mathcal{M}_n)$ and $D_n$ as well:
in particular $i^m_n : T_p(\mathcal{M}_m) \hookrightarrow T_p(\mathcal{M}_n)$ is injective and satisfies
$T_p(\mathcal{M}_n)^{\Gamma_m^n} = i^m_n \left( T_p(\mathcal{M}_m) \right)$
(see \cite[Theorem 3.1]{Greither Popescu Galois}) and the composition $N^n_m \circ i^m_n$ is multiplication 
by $[F_n : F_m]$ (on both modules).
All these maps are compatible with the exact sequence (\ref{successione esatta Greither Popescu}), i.e. 
the following diagram is commutative for every pair of indices $n>m$
\[ \xymatrix{ 0 \ar[r] & T_p(F_n) \ar[rr] \ar@/_10pt/[dd]_{N^n_m} & & T_p(\mathcal{M}_n) \ar[rr] \ar@/_10pt/[dd]_{N^n_m} & &
D_n \ar[r] \ar@/_10pt/[dd]_{N^n_m} & 0 \\
\  \\
0 \ar[r] & T_p(F_m) \ar[rr] \ar@/_10pt/[uu]_{i^m_n} & & T_p(\mathcal{M}_m) \ar[rr] \ar@/_10pt/[uu]_{i^m_n} & &
D_m \ar[r] \ar@/_10pt/[uu]_{i^m_n} & 0\,. } \]

\begin{lemma}\label{LemNormTM}
    The map $N^n_m : T_p(\mathcal{M}_n)\rightarrow T_p(\mathcal{M}_m)$ is surjective,
    its kernel is $I_{\Gamma^n_m}T_p(\mathcal{M}_n)$.
\end{lemma}

\begin{proof}
    By \cite[Theorem 3.9 part (2)]{Greither Popescu Galois} $T_p(\mathcal{M}_n)$
    is a free $\Zp[\Gamma^n_m]$-module (because $\Gamma^n_m$ is a $p$-group), hence cohomologically trivial, i.e.
    \begin{displaymath}
        \widehat{H}^i\left(\Gamma^n_m, T_p(\mathcal{M}_n) \right) = 0 \quad \text{for every integer } i,
    \end{displaymath}
    where the $\widehat{H}^i(\Gamma^n_m, \bullet)$ are Tate cohomology groups.
    
    \noindent Specializing the previous equality at $i=0$ we obtain 
    \begin{displaymath}
        N^n_m \left( T_p(\mathcal{M}_n) \right) = T_p(\mathcal{M}_n)^{\Gamma_m^n}
            = i^m_n \left( T_p(\mathcal{M}_m) \right)\simeq T_p(\mathcal{M}_m).
    \end{displaymath}
    In a similar way we obtain the second part of the lemma by specializing at $i=-1$.
\end{proof}

We can now prove an analogous result for the modules $T_p(F_n)$.

\begin{prop}\label{proposizione su norm di T_p}
    Let $\chi$ be a character of type $1$ or $2$. Then the norm $N^n_m:T_p(F_n)(\chi)\rightarrow T_p(F_m)(\chi)$ is surjective and
    its kernel is $I_{\Gamma^n_m}T_p(F_n)(\chi)$.
\end{prop}

\begin{proof}
    Consider the $\chi$-part of the exact sequence \eqref{successione esatta Greither Popescu}, i.e.
    \begin{equation}\label{successione esatta Greither Popescu chi-parte}
    \xymatrix{ 0\ar[r] & T_p(F_n)(\chi) \ar[r] & T_p(\mathcal{M}_n)(\chi) \ar[r] & D_n(\chi) \ar[r] & 0 \,.}
    \end{equation} 
    If $\chi$ is a character of type $1$, then $T_p(F_n)(\chi) \simeq T_p(\mathcal{M}_n)(\chi)$ and this is a 
    restatement of Lemma \ref{LemNormTM}. 
    
    \noindent If $\chi$ is of type $2$, then $D_n(\chi) \simeq W[\Gamma_n]$ so both $D_n(\chi)$ and
    $T_p(\mathcal{M}_n)(\chi)$ are cohomologically trivial. Hence 
    \begin{equation}\label{T_p coomologicamente banale}
        \widehat{H}^i\left(\Gamma_m^n, T_p(F_n)(\chi) \right) = 0 \quad \text{for every integer}\ i \,
    \end{equation}
    and, specializing at $i=-1$, we obtain the statement for the kernel.
    
    \noindent To prove surjectivity we take $\Gamma_m^n$-invariants in
    \eqref{successione esatta Greither Popescu chi-parte} to obtain
 \[ \xymatrix{ 0\ar[r] & T_p(F_n)(\chi)^{\Gamma^n_m} \ar[r] & T_p(\mathcal{M}_m)(\chi) \ar[r] & 
D_n(\chi)^{\Gamma^n_m} \ar[r] & \widehat{H}^1\left(\Gamma_m^n , T_p(F_n)(\chi) \right)  \,.} \]
    Note that $D_n(\chi)^{\Gamma^n_m} \simeq W[\Gamma_n]^{\Gamma^n_m} \simeq W[\Gamma_m]\simeq D_m(\chi)$,
    so, comparing with \eqref{successione esatta Greither Popescu chi-parte}
    with $m$ in place of $n$, we get $T_p(F_n)(\chi)^{\Gamma^n_m} \simeq T_p(F_m)(\chi)$. Now
    specializing \eqref{T_p coomologicamente banale} at $i=0$ we obtain
    \begin{displaymath}
        N^n_m \left( T_p(F_n)(\chi)\right) = T_p(F_n)(\chi)^{\Gamma^n_m} \simeq T_p(F_m)(\chi) \,.\qedhere
    \end{displaymath}
\end{proof}

We are now ready to prove the main theorems of this section.

\begin{teo}\label{ThmTateFinGen}
    Let $\chi$ be a character of type $1$ or $2$. Then $T_p(\mathcal{F})(\chi)$ is a finitely generated
    $\Lambda_\F$-module. Moreover, if $\displaystyle{\lim_\leftarrow \Theta_n^{\sharp}(\gamma^{-1},\chi)\neq 0}$, 
    then $T_p(\mathcal{F})(\chi)$ is torsion.
\end{teo}

\begin{proof}
    Fix $m$, and let $\mathfrak{I}_m$ be the augumentation ideal of $\Gal(\mathcal{F} / F_m)$ in
    $\Lambda:=W\llbracket \Gamma_\infty \rrbracket$. Let $\widetilde{\mathfrak{I}}_m =
    \Lambda_\F \otimes_{W\llbracket \Gamma_\infty \rrbracket} \mathfrak{I}_m$ be the corresponding ideal in 
    $\Lambda_\F=\Lambda\llbracket G_\F\rrbracket$, and recall that 
    $\displaystyle{ \mathfrak{I}_m = \lim_{\substack{\leftarrow \\ n}} I_{\Gamma^n_m}}$.\\
    By Proposition \ref{proposizione su norm di T_p} we have 
    \begin{displaymath}
      T_p(F_m)(\chi) = N^n_m \left(T_p(F_n)(\chi)\right) \simeq T_p(F_n)(\chi) / I_{\Gamma^n_m}T_p(F_n)(\chi)\,.
    \end{displaymath}
    This holds for every $n >m$ and so 
    \begin{displaymath}
        T_p(F_m)(\chi) \simeq T_p(\mathcal{F})(\chi) / \widetilde{\mathfrak{I}}_m T_p(\mathcal{F})(\chi) \,.
    \end{displaymath}
    The module on the left is finitely generated over $W[\Gamma_m]\llbracket G_\F \rrbracket = \Lambda_\F / \widetilde{\mathfrak{I}}_m$
    and, since the ideals $\widetilde{\mathfrak{I}}_m$ form an open filtration of the profinite
    algebra $\Lambda_\F$, we can apply the generalized Nakayama Lemma of \cite{BH} and obtain that
    $T_p(\mathcal{F})(\chi)$ is a finitely generated $\Lambda_\F$-module. \\
    Now we define the element $\Theta_\infty^{\sharp}(\gamma^{-1},\chi) \in \Lambda_\F$ as
    \begin{displaymath}
        \Theta_\infty^{\sharp}(\gamma^{-1},\chi) =
        \left\{ \begin{array}{ll}
            \Theta_\infty(\gamma^{-1},\chi)                         &   \ \mbox{if } \chi \mbox{ is of
                                                                                    type } 1 \,, \\
                \ & \\                                                                    
                \displaystyle{\frac{\Theta_\infty(\gamma^{-1},\chi)}{1-\gamma^{-1}} }   &   \ \mbox{if } \chi \mbox{ is of
                                                                                    type } 2 \,,
        \end{array} \right.
    \end{displaymath}
    i.e. the inverse limit of the generators of ${\rm Fitt}_{W[\Gamma_m]\llbracket G_\F \rrbracket}(T_p(F_m)(\chi))$,
    by Theorem \ref{Teorema fitting per caratteri di tipo 1 o 2}. 
	Clearly $\Theta_\infty^{\sharp}(\gamma^{-1},\chi)T_p(\mathcal{F})(\chi) = 0$ and this implies the statement on torsion.
\end{proof}

\noindent Now we know that the Fitting ideal of $T_p(\mathcal{F})(\chi)$ is well defined and we proceed to compute a generator via a
limit process.

\begin{teo}\label{TeoFittIdealTateIwasawaModule}
    Let $\chi$ be a character of type $1$ or $2$. Then
    \begin{displaymath}
        {\rm Fitt}_{\Lambda_\F} \left(T_p(\mathcal{F})(\chi) \right) =
            \left( \Theta_\infty^{\sharp}(\gamma^{-1},\chi)\right) \,.
    \end{displaymath}
\end{teo}

\begin{proof}
	The equalities
	\[
		\left( \Theta_\infty^{\sharp}(\gamma^{-1},\chi)\right) = \lim_\leftarrow
		\left( \Theta_n^{\sharp}(\gamma^{-1},\chi)\right) = \lim_\leftarrow
		{\rm Fitt}_{W[\Gamma_n]\llbracket G_\F \rrbracket} \left( T_p(F_n)(\chi) \right)
	\]
	reduce the proof to showing 
	\[
		{\rm Fitt}_{\Lambda_\F} \left(T_p(\mathcal{F})(\chi) \right) = \lim_\leftarrow 
		{\rm Fitt}_{W[\Gamma_n]\llbracket G_\F \rrbracket} \left( T_p(F_n)(\chi) \right)\,.
	\]
    Let $N^\infty_m : T_p(\mathcal{F})(\chi) \twoheadrightarrow T_p(F_m)(\chi)$ be induced by the projection modulo 
    $\widetilde{\mathfrak{I}}_mT_p(\mathcal{F})(\chi)$. These maps are obviously compatible with the norm maps, i.e. 
    $N^\infty_m =N^n_m \circ N^\infty_n$ for any $n>m$. Let $t_1, \dots ,t_r$ be $\Lambda_\F$-generators of $T_p(\mathcal{F})(\chi)$
    and write $K_\infty$ for the kernel of the surjective map $\Lambda_{\mathbb{F}}^{\oplus r} \twoheadrightarrow T_p(\mathcal{F})(\chi)$
    sending the $i$-th element of the canonical basis to $t_i$. 
    We have an exact sequence 
\begin{equation}\label{EqExSeqInfty}
\xymatrix{ 0 \ar[r] & K_\infty \ar[r] & \Lambda_\F^{\oplus r} \ar[r] &
T_p(\mathcal{F})(\chi) \ar[r] & 0 \,.} \end{equation}
    Since $N^\infty_n(t_1), \dots , N^\infty_n(t_r)$ generate $T_p(F_n)(\chi)$ over $W[\Gamma_n]\llbracket G_\F \rrbracket$, we can construct
    similar exact sequences for every integer $n$, i.e.
\[ \xymatrix{ 0 \ar[r] & K_n \ar[r] & W[\Gamma_n]\llbracket G_\F \rrbracket^{\oplus r} \ar[r] &
T_p(F_n)(\chi) \ar[r] & 0 } \]
   (where now the $i$-th element of the canonical basis maps to $N^\infty_n(t_i)$). They all fit into the commutative diagrams
\[ \xymatrix{ 0 \ar[r] & K_n \ar[rr] \ar[dd]^{k^n_m} & & W[\Gamma_n]\llbracket G_\F \rrbracket^{\oplus r} \ar[rr] \ar[dd]^{\pi^n_m} 
& & T_p(F_n)(\chi) \ar[r] \ar[dd]^{N^n_m} & 0 \\
\  \\
0 \ar[r] & K_m \ar[rr] & & W[\Gamma_m]\llbracket G_\F \rrbracket^{\oplus r} \ar[rr]  
& & T_p(F_m)(\chi) \ar[r] & 0 \,, }\]
    where $k^n_m:=(\pi^n_m)_{|K_n}$.

\noindent One has $\Ker(\pi^n_m)=(I_{\Gamma^n_m}W[\Gamma_n]\llbracket G_\F \rrbracket)^{\oplus r}$ and 
$\Ker(N^n_m)=I_{\Gamma^n_m}T_p(F_n)(\chi)$ (by Proposition \ref{proposizione su norm di T_p}), so the map between them 
    is surjective. Since $\pi^n_m$ is surjective, the snake lemma implies that $k^n_m$ is surjective as well and the diagram 
    above satisfies the Mittag-Leffler condition which allows us to
    take the inverse limit. Comparing this limit with \eqref{EqExSeqInfty}, we obtain 
    $\displaystyle{ K_\infty = \lim_\leftarrow K_n}$.
    
	\noindent Let $\mathscr{M}_r(K_m)$ be the set of $r\times r$ matrices whose entries are
	in $W[\Gamma_m]\llbracket G_\F \rrbracket$ and such that each row, seen as a vector in 
	$W[\Gamma_m]\llbracket G_\F \rrbracket^{\oplus r}$, is in $K_m$; we still denote by $k^n_m$ the natural
	extension of the map $K_n \rightarrow K_m$ to $\mathscr{M}_r(K_n) \rightarrow
	\mathscr{M}_r(K_m)$, which is surjective as well. By definition ${\rm Fitt}_{W[\Gamma_m]\llbracket G_\F \rrbracket} 
	\left( T_p(F_m)(\chi) \right)$ is the ideal generated by $\{ \det(M_m)\,:\,M_m \in \mathscr{M}_r(K_m)\}$.
	The commutativity of the previous diagram yields for each $M_n\in \mathscr{M}_r(K_n)$
	\[
		\pi^n_m \left( \det(M_n) \right) = \det \left( k^n_m (M_n) \right).
	\]
	Extending this construction to the infinite level (with analogous notations), we obtain
	$\pi^\infty_m \left( \det (M_\infty) \right) \in 
	{\rm Fitt}_{W[\Gamma_m]\llbracket G_\F \rrbracket} \left( T_p(F_m)(\chi) \right)$ for any $M_\infty \in \mathscr{M}_r(K_\infty)$.
	Hence, for any $m$, $\pi^\infty_m \left(  {\rm Fitt}_{\Lambda_\F} \left(T_p(\mathcal{F})(\chi) \right)\right)
	\subseteq {\rm Fitt}_{W[\Gamma_m]\llbracket G_\F \rrbracket} \left( T_p(F_m)(\chi) \right)$, and 
	\[
		{\rm Fitt}_{\Lambda_\F} \left(T_p(\mathcal{F})(\chi) \right) \subseteq \lim_\leftarrow 
		{\rm Fitt}_{W[\Gamma_m]\llbracket G_\F \rrbracket} \left( T_p(F_m)(\chi) \right).		
	\]
	The other inclusion needs a little bit more work, basically
	we follow the arguments of \cite[Theorem 2.1]{Greither Kurihara}.
	Each element of ${\rm Fitt}_{W[\Gamma_m]\llbracket G_\F \rrbracket} \left( T_p(F_m)(\chi) \right)$
	can be written as a linear combination
	\begin{equation}\label{EqCombLinDet}
		x_m = \sum_{i=1}^s \lambda_i \, \det (M_m^{(i)})
	\end{equation}
	with $\lambda_i \in W[\Gamma_m]\llbracket G_\F \rrbracket$ and $M_m^{(i)} \in \mathscr{M}_r(K_m)$;
	multiplying the first row of the $M_m^{(i)}$ by $\lambda_i$, we get matrices $\overline{M}_m^{(i)}$ such that
	\[	x_m = \sum_{i=1}^s \det (\overline{M}_m^{(i)}) \,,\]
 	i.e. we can assume that all coefficients in \eqref{EqCombLinDet} are $1$. Since the number of elements needed to generate
	$T_p(F_m)(\chi)$ (and $T_p(\mathcal{F})(\chi)$) can be chosen independently from $m$, one has that
	$s$ can be chosen independently of $m$ as well.
	
	\noindent Now put $\mathcal{B}_m := \mathscr{M}_r(K_m)^{\oplus s}$ with the induced topology
	and define the non-linear operator $\phi_m : \mathcal{B}_m \rightarrow W[\Gamma_m]\llbracket G_\F \rrbracket$, by
	$\phi_m \left( M_m^{(1)}, \dots, M_m^{(s)} \right) = \sum_i \det (M_m^{(i)})$
	(analogous definition for $\phi_\infty$). This operator is continuous and its image is  
	${\rm Fitt}_{W[\Gamma_m]\llbracket G_\F \rrbracket}T_p(F_m)(\chi)$. We extend the map $k^n_m$ from $\mathscr{M}_r(K_n)$
	to $\mathcal{B}_n$ and get the commutative diagram
	\[
		\xymatrix{\mathcal{B}_n \ar[r]^(.35){\phi_n} \ar[d]_{k^n_m} & W[\Gamma_n]\llbracket G_\F \rrbracket \ar[d]^{\pi^n_m} \\
				  \mathcal{B}_m \ar[r]_(.35){\phi_m} 		  & W[\Gamma_m]\llbracket G_\F \rrbracket \, .}
	\]
	Now we take a sequence $\displaystyle{ (x_m)_{m \in \N} \in \lim_\leftarrow \, 
		{\rm Fitt}_{W[\Gamma_m]\llbracket G_\F \rrbracket} \left( T_p(F_m)(\chi) \right) }$,
	and look for an element $b_\infty \in \mathcal{B}_\infty$ such that $(x_m)_{m \in \N}=\phi_\infty (b_\infty)\in 
	{\rm Fitt}_{\Lambda_\F} \left(T_p(\mathcal{F})(\chi) \right)$. \\
	For any $m$ put $\Upsilon_m := \phi_m^{-1}(x_m)$, then $\Upsilon_m$ is closed and, since $W[\Gamma_m]\llbracket G_\F \rrbracket$ is compact,
	$\Upsilon_m$ is compact as well. For each $\upsilon_n \in \Upsilon_n$ we have that
	\[
		\phi_m \left( k^n_m(\upsilon_n) \right) = \pi^n_m \left( \phi_n(\upsilon_n) \right) = \pi^n_m(x_n) = x_m\,,
	\]
	thus $k^n_m(\Upsilon_n)\subseteq \Upsilon_m$ and we define 
	\[
		\overline{\Upsilon}_m = \bigcap_{n > m} k^n_m (\Upsilon_n) \subseteq \Upsilon_m.
	\]
	Since $k^n_m$ is continuous, $\overline{\Upsilon}_m$ is compact and not empty, moreover one easily shows that
	$k^n_m(\overline{\Upsilon}_n)\subseteq \overline{\Upsilon}_m$ and we are going to prove equality between them.
	Let $\overline{\upsilon}_m \in \overline{\Upsilon}_m$ so that, for any $n>m$, there exists $\upsilon_n \in \Upsilon_n$ with $k^n_m(\upsilon_n) = \overline{\upsilon}_m$. 
	Now fix $\ell>0$ and, for $n>m+\ell$, consider $k^n_{m+\ell}(\upsilon_n) \in \Upsilon_{m+\ell}$ as a sequence in $n$. Since 
	$\Upsilon_{m+\ell}$ is compact, there exists a convergent subsequence $\upsilon_{n_j}$ whose limit we call $\overline{\upsilon}_{m+\ell}$.
Then, for any $n>m+\ell$, 	
\[ \overline{\upsilon}_{m+\ell} = \lim_{j \rightarrow \infty } k^{n_j}_{m+\ell}(\upsilon_{n_j}) = 
\lim_{\substack{j \rightarrow \infty \\ n_j\geqslant n } }	( k^n_{m+\ell}\circ k^{n_j}_n)(\upsilon_{n_j})
\in k^n_{m+\ell}(\Upsilon_n) \,, \]
i.e. $\overline{\upsilon}_{m+\ell}$ is in $\overline{\Upsilon}_{m+\ell}$. Obviously $k^{m+\ell}_m(\overline{\upsilon}_{m+\ell})=\overline{\upsilon}_m$, 
so the map $k^{m+\ell}_m$ is surjective and we have constructed a coherent sequence $b_\infty:=\overline{\upsilon}_\infty=(\overline{\upsilon}_m)_{m\in \N}
\in \mathcal{B}_\infty$. Since $\phi_m(\overline{\upsilon}_m) = x_m$ for each integer $m$, we have that 
$\phi_\infty (b_\infty) = (x_m)_{m \in \N}$.
\end{proof}

\subsection{Fitting ideals for the class groups: finite level}\label{SecClassGroup}
Now we move to our primary interest: the $p$-part $\mathcal{C}\ell^0(F_n)$ of the class groups 
of degree zero divisors of the field $F_n$ which, from now on, we shall denote by $C_n$ and
which is naturally a finitely generated torsion $W[\Gamma_n]$-module. To compute its Fitting ideal
we shall use specializations of ${\rm Fitt}_{W[\Gamma_n]\llbracket G_\F\rrbracket}(T_p(F_n)(\chi))$ 
as suggested by the following lemmas (for the first see e.g. \cite[Lemma 4.6]{ABBL}, the second is a well known property
of Fitting ideals). 

\begin{lemma}\label{LemComputoFittClassGr1}
	There is an isomorphism of $\Zp[G_n]$-modules
	\[ T_p(F_n)_{G_\F}:=T_p(F_n) / (1-\gamma^{-1})T_p(F_n)\simeq C_n.\]
\end{lemma}

\begin{lemma}\label{LemComputoFittClassGr2}
	Let $M$ be a finitely generated torsion module over $R$. Let
	$I$ be any nontrivial ideal of $R$ and consider the projection $\pi_I : R \twoheadrightarrow R/ I$. Then
	\[ {\rm Fitt}_{R/I} \left( M/IM \right) = \pi_I \left( {\rm Fitt}_R (M) \right). \]
\end{lemma}

\noindent Let $\pi : W[\Gamma_n] \llbracket G_\F \rrbracket \rightarrow W[\Gamma_n]\simeq W[\Gamma_n] \llbracket G_\F \rrbracket / I_{G_\F}$ 
be the canonical projection sending $\gamma$ to $1$. Combining Lemmas \ref{LemComputoFittClassGr1} and 
\ref{LemComputoFittClassGr2}, and the computations of Theorem \ref{Teorema fitting per caratteri di tipo 1 o 2} we obtain

\begin{teo}\label{Teorema fitting Class Gr per caratteri di tipo 1 o 2} 
    Let $\chi \in \widehat{G_0}$ be a character of type 1 or 2. Then
    \[ {\rm Fitt}_{W[\Gamma_n]} \left( C_n(\chi)\right) = \left( \Theta_n^\sharp (1,\chi) \right), \]
    where   
    \begin{displaymath}
        \Theta_n^\sharp (1,\chi) =
            \left\{ \begin{array}{ll}
                \Theta_n(1,\chi)                        &   \ \mbox{if } \chi \mbox{ is of
                                                                                    type } 1\,, \\
                 \ & \\                                                                   
                \displaystyle{\bigg(\frac{\Theta_n(\gamma^{-1},\chi)}{1-\gamma^{-1}}\bigg)_{\bigr| \gamma=1} }  &   \ \mbox{if } \chi \mbox{ is of
                                                                                    type } 2\,.
            \end{array} \right.
    \end{displaymath}
\end{teo}

\begin{proof}
Just specialize Theorem \ref{Teorema fitting per caratteri di tipo 1 o 2} to $\gamma=1$ recalling the convergence properties of the
Stickelberger series.
\end{proof}

\subsection{Fitting ideals for the class groups: infinite level and the Main Conjecture}
Now we approach the $W\llbracket \Gamma\rrbracket$-module $C_\infty := \displaystyle{ \lim_\leftarrow C_n }$,
where the limit is with respect to the norm maps $N^n_m : C_n \rightarrow C_m$.
We shall also consider maps $i^m_n : C_m \rightarrow C_n$ induced by the embeddings $i^m_n : {\rm Div}(F_m) \rightarrow {\rm Div}(F_n)$. 
We recall that for any $D=\sum_\nu n_\nu \nu\in {\rm Div}(F_m)$, we have 
$i^m_n (D) := \displaystyle{\sum_\nu n_\nu \sum_{w | \nu} e(w | \nu) w}$, where $ e(w | \nu)$ is the ramification index of
$w$ over $\nu$. In particular, since $	\deg \left( i^m_n (D) \right) =[F_n : F_m]\cdot \deg(D)$, the image of a degree zero divisor 
still has degree zero, moreover $i^m_n({\rm Div}(F_m))$ is $\Gamma^n_m$-invariant. 

\noindent The following proposition gives us information on injectivity and surjectivity of $N$ and $i$.

\begin{prop}\label{PropNormSurjImmerInject} Let $F_0 \subseteq K \subset E \subset \mathcal{F}$ with
	$[E:F_0]$ finite. Then
	\begin{enumerate}
	\item[(a)] the norm map $N^E_K : \mathcal{C}\ell^0(E) \rightarrow \mathcal{C}\ell^0(K)$ is surjective;
	\item[(b)] the map $i^K_E : \mathcal{C}\ell^0(K) \rightarrow \mathcal{C}\ell^0(E)$ is injective.
	\end{enumerate}
\end{prop}
\begin{proof}
	We recall that $d_\infty=1$ and that $\infty$ is totally split in $\mathcal{F}/F_0$, hence all primes in $E$ and $K$ dividing $\infty$
	have degree 1 as well.
	 
\noindent 	(a)
This is just an application of class field theory for function fields, see
	e.g. \cite[Lemma 5.4 part (3)]{ABBL}.
	
\noindent (b)
	Let $G:= {\rm Gal} (E/K)$ and, for any field $L$, let $\mathcal{P}_L$ be the 
	principal divisors of $L$. Taking the $G$-cohomology in the exact sequence
	\begin{equation}\label{SuccEsattaDivPrincip}
		\xymatrix{ 0 \ar[r] & \Fq^\times \ar[r] & E^\times \ar[r] & \mathcal{P}_E \ar[r] &  0 \,,}
	\end{equation}
	we get
	\[ 		\xymatrix{ 0 \ar[r] & \Fq^\times \ar[r] & K^\times \ar[r] & \mathcal{P}_E^G \ar[r] &  0 }\]
	and 	$H^1(G,\mathcal{P}_E) = 0$
	(because of Hilbert 90 and $(|G|,\Fq^\times)=1$).
	Comparing this with the analogue of \eqref{SuccEsattaDivPrincip} for $K$
	we have that $\mathcal{P}_E^G = \mathcal{P}_K$. Taking the $G$-cohomology in
	\[
		\xymatrix{ 0 \ar[r] & \mathcal{P}_E \ar[r] & {\rm Div}^0(E) \ar[r] & \mathcal{C}\ell^0(E) \ar[r] &  0 ,}
	\]
	we obtain
	\[
		\xymatrix{ 0 \ar[r] & \mathcal{P}_E^G = \mathcal{P}_K \ar[r] & {\rm Div}^0(E)^G \ar[r] & \mathcal{C}\ell^0(E)^G \ar[r] & 
				H^1(G, \mathcal{P}_E) ,}\]
	which fits into the following commutative diagram 
	\begin{equation}\label{DiagrammaG-invClassGr}
	\xymatrix{ 0 \ar[r] & \mathcal{P}_K \ar[r]\ar@{=}[d] & {\rm Div}^0(K) \ar[r]\ar@{^(->}[d] & \mathcal{C}\ell^0(K) \ar[r]\ar[d]^{i^K_E} & 
				0 \\
			0 \ar[r] & \mathcal{P}_K \ar[r] & {\rm Div}^0(E)^G \ar[r] & \mathcal{C}\ell^0(E)^G \ar[r] & 0.}
	\end{equation}	
	Applying the snake lemma we obtain the thesis. \qedhere
\end{proof}

\noindent From diagram \eqref{DiagrammaG-invClassGr} we also deduce that
\begin{equation}\label{EqIsomQuotClassGr}
	\mathcal{C}\ell^0(E)^G / i^K_E \left( \mathcal{C}\ell^0(K) \right) \simeq {\rm Div}^0(E)^G / i^K_E \left( {\rm Div}^0(K)\right)\,.
\end{equation} 

\noindent To perform a limit and move to the infinite level we still have to deal with the kernel of the norm map 
$N^E_K(\chi) : \mathcal{C}\ell^0(E)(\chi) \rightarrow \mathcal{C}\ell^0(K)(\chi)$ for characters of type 1 or 2.
 
\begin{lemma}
	 Let $F_0 \subseteq K \subset E \subset \mathcal{F}$ with $[E:F_0]$ finite and let $G:= {\rm Gal} (E/K)$. 
	 Assume $|G|=p$, then the group $\Delta := \Gal(F_0 / H_A)$ acts trivially on
	 $\mathcal{C}\ell^0(E)^G / i^K_E \left( \mathcal{C}\ell^0(K) \right)$.
\end{lemma}

\begin{proof}
	The composition of two natural maps
	\[
	{\rm Div}^0(E)^G \hookrightarrow {\rm Div}(E)^G \twoheadrightarrow {\rm Div}(E)^G / i^K_E \left( {\rm Div}(K) \right)\,,
	\]
has kernel ${\rm Div}^0(E)^G \cap i^K_E \left( {\rm Div}(K) \right) = 
	i^K_E \left( {\rm Div}^0(K) \right)$ and induces an injection
	\[
		{\rm Div}^0(E)^G / i^K_E \left( {\rm Div}^0(K)\right) \hookrightarrow 
		{\rm Div}(E)^G / i^K_E \left( {\rm Div}(K) \right).
	\]
	Thus it is enough to show that $\Delta$ acts trivially on ${\rm Div}(E)^G / i^K_E \left( {\rm Div}(K) \right)$. 
	\noindent	Let $\mathfrak{p}_1, \dots, \mathfrak{p}_s$ (resp. $\mathfrak{P}_1, \dots, \mathfrak{P}_s$) be the set of primes of 
	$K$ (resp. $E$) lying above $\mathfrak{p}$. The extension $E/K$ is totally ramified at $\mathfrak{p}$ and we can assume 
	that $\mathfrak{P}_j$ is the unique prime of $E$ lying above $\mathfrak{p}_j$, i.e. $i^K_E (\mathfrak{p}_j) = p \mathfrak{P}_j$.
	Moreover the only extension where the prime $\mathfrak{p}$ may split is $H_A /F$, so $s$ divides $h^0(F)$, i.e. is coprime
	with $p$. 
	\noindent	We can write ${\rm Div}(K) = \bigoplus_\nu \Z \nu$ (where $\nu$ runs through all the primes of $K$)
	and ${\rm Div}(E) = \bigoplus_\nu H_\nu$, with $H_\nu = \bigoplus_{w \mid \nu} \Z w$. Now for the ramified primes we have
	$H_{\mathfrak{p_j}} = \Z \mathfrak{P}_j = H_{\mathfrak{p}_j}^G$ while, for the unramified ones, 
	if we let $G_\nu$ be the decomposition	group of $\nu$ in $G$, we have $H_\nu = \Z [G / G_\nu]w$, so that 
	$H_\nu^G = i^K_E \left(\Z \nu \right)$. Therefore
	\[
		{\rm Div}(E)^G = \bigoplus_{j=1}^s \Z \mathfrak{P}_j \oplus \bigoplus_{\nu \nmid \mathfrak{p}} 
		i^K_E \left(\Z \nu \right) \,,
	\]
	\[
		i^K_E \left( {\rm Div}(K) \right) = \bigoplus_{j=1}^s p\Z \mathfrak{P}_j \oplus \bigoplus_{\nu \nmid \mathfrak{p}} 
		i^K_E \left(\Z \nu \right)
	\]
	and finally
	\[
		{\rm Div}(E)^G / i^K_E \left( {\rm Div}(K) \right) = \bigoplus_{j=1}^s \left( \Z / p\right) \mathfrak{P}_j \,.
	\]
	Note that, for any set of integers $\alpha_1 , \dots , \alpha_s$ coprime with $p$, the classes  
	$\alpha_j \mathfrak{P}_j$ with $j=1,\dots,s$ still generate ${\rm Div}(E)^G / i^K_E \left( {\rm Div}(K) \right)$.
	
	\noindent Now consider the subfield $E^\Delta$ of $E$ (resp. $K^\Delta$ of $K$) fixed by $\Delta$. 
	Since $|\Delta|$ is prime with $p$, there is a canonical isomorphism $G^\Delta:= \Gal(E^\Delta / K^\Delta) \simeq G$ and,
	since $F_0 / H_A$ is totally ramified at $\mathfrak{p}$, we still have exactly $s$ primes in $E^\Delta$ (resp. $K^\Delta$) 
	above $\mathfrak{p}$: let $\mathfrak{P}_j^\Delta$ (resp. $\mathfrak{p}_j^\Delta$) be those primes and, as above,
	assume $i^{K^\Delta}_ {E^\Delta} (\mathfrak{p}_j^\Delta) = p \mathfrak{P}_j^\Delta$. 
	With the same argument, we can prove
	\[
		{\rm Div}(E^\Delta)^{G^\Delta} / i^{K^\Delta}_{E^\Delta} \left( {\rm Div}(K^\Delta) \right) = 
		\bigoplus_{j=1}^s \left( \Z / p\right) \mathfrak{P}_j^\Delta \,.
	\]
	To conclude note that $i_E^{E^\Delta} (\mathfrak{P}^\Delta_j)
	= |\Delta | \mathfrak{P}_j$; since $|\Delta |$ is coprime with $p$, these classes generate  
	${\rm Div}(E)^G / i^K_E \left( {\rm Div}(K) \right)$
	and clearly the action of $\Delta$ on them is trivial.
\end{proof}

\noindent We use the previous lemma to prove

\begin{prop}\label{PropKernelNormClGr}
	 Let $F_0 \subseteq K \subset E \subset \mathcal{F}$ with $[E:F_0]$ finite and let $G:= {\rm Gal} (E/K)$. 
	 If $\chi\in \widehat{G_0}$ is of type $1$ or $2$, then 
	 \[ \Ker\left( N^E_K(\chi) \right) = I_G \mathcal{C}\ell^0(E)(\chi)\,, \]
	 where $I_G$ denotes the augumentation of $G$.
\end{prop}

\begin{proof}
	We proceed by induction on $|G|$, starting with the case $|G|=p$ (for $|G|=1$ there is nothing to prove).
	Since $\chi$ is not of type 3, it may be seen as a nontrivial character of 
	$\Delta = \Gal(F_0 / H_A)$ and, by the previous lemma, 
	\[ \left(\mathcal{C}\ell^0(E)^G / i^K_E \left( \mathcal{C}\ell^0(K) \right)\right)(\chi) = 0,\quad
	{\rm i.e.}\quad \mathcal{C}\ell^0(E)^G (\chi) = i^K_E \left( \mathcal{C}\ell^0(K) \right)(\chi) .\] 
	Let $g$ be a generator of $G$, then $I_G \mathcal{C}\ell^0(E)(\chi) = (1-g)\mathcal{C}\ell^0(E)(\chi)$ and we also recall that 
	the cyclicity of $G$ yields $I_G \mathcal{C}\ell^0(E)(\chi)\simeq \mathcal{C}\ell^0(E)/\mathcal{C}\ell^0(E)^G$. 
	We have two exact sequences
	\[\xymatrix{ 0 \ar[r]& \Ker\left( N^E_K(\chi) \right) \ar[r] & \mathcal{C}\ell^0(E)(\chi) \ar[r]^{N^E_K(\chi)\ } &
			\mathcal{C}\ell^0(K)(\chi) \ar[r] & 0,} 
	\]	
	(exact by Proposition \ref{PropNormSurjImmerInject}, part (a)), and
	\[\xymatrix{ 0 \ar[r]& \mathcal{C}\ell^0(K)(\chi) \ar[r]^{i^K_E} & \mathcal{C}\ell^0(E)(\chi) \ar[r]^{\!\!\!1-g\ } &
			I_G \mathcal{C}\ell^0(E)(\chi) \ar[r] & 0} 
	\]	
	(exact by what we noted above). Cardinalities yield $|\Ker\left( N^E_K(\chi) \right)| = |I_G \mathcal{C}\ell^0(E)(\chi)|$,
	and, since $I_G \mathcal{C}\ell^0(E)(\chi) \subseteq \Ker\left( N^E_K(\chi) \right)$, we have equality of the two groups. 
	
	\noindent  For the inductive step assume $|G|=p^l > p$ and take an intermediate field 
	$K \subsetneqq E' \subsetneqq E$; put $G_1=\Gal(E/E')$ and $G_2=\Gal(E'/K)$, so that $G_1$ and 
	$G_2$ have cardinality strictly smaller than $p^l$. The inductive hypothesis yields
	\[ \Ker\left( N^E_{E'}(\chi) \right) = I_{G_1} \mathcal{C}\ell^0(E)(\chi) \]
	and
	\[ \Ker\left( N^{E'}_K(\chi) \right) = I_{G_2} \mathcal{C}\ell^0(E')(\chi). \]
	By Proposition \ref{PropNormSurjImmerInject} part (a), the norm $N^E_{E'} : \mathcal{C}\ell^0(E)(\chi) \rightarrow 
	\mathcal{C}\ell^0(E')(\chi)$ is surjective and so
	\[	N^E_{E'}\left( I_G \mathcal{C}\ell^0(E)(\chi) \right) = I_{G_2} \mathcal{C}\ell^0(E')(\chi). \]
	Let $x \in \Ker\left( N^E_K(\chi) \right)$, since $N^E_K = N^{E'}_K \circ N^E_{E'}$, 
	we have that $N^E_{E'}(x)\in \Ker\left( N^{E'}_K(\chi) \right) = I_{G_2} \mathcal{C}\ell^0(E')(\chi)$, so
	there exists $\alpha \in I_G$ and $y \in \mathcal{C}\ell^0(E)(\chi)$ such that $N^E_{E'}(x) = N^E_{E'}(\alpha y)$.
	Therefore $x - \alpha y \in \Ker\left( N^E_{E'}(\chi) \right)$, which yields
	\[ x \in \Ker\left( N^E_{E'}(\chi) \right) + I_G \mathcal{C}\ell^0(E)(\chi) = I_G \mathcal{C}\ell^0(E)(\chi), \]
	since $\Ker\left( N^E_{E'}(\chi) \right) = I_{G_1} \mathcal{C}\ell^0(E)(\chi) \subseteq I_G \mathcal{C}\ell^0(E)(\chi)$.
	We have proved that 
	\[ \Ker\left( N^E_K(\chi) \right) \subseteq I_G \mathcal{C}\ell^0(E)(\chi) \]
	and the other inclusion is trivial.
\end{proof}

We can finally prove the main theorems on $C_\infty$.

\begin{teo}\label{ThmClassGrFinGen}
    Let $\chi\in\widehat{G_0}$ be a character of type $1$ or $2$. Then $C_\infty(\chi)$ is a finitely generated
    $\Lambda:=W\llbracket \Gamma_\infty\rrbracket$-module. Moreover, let
    \begin{displaymath}
        \Theta_\infty^{\sharp}(1,\chi) =
        \left\{ \begin{array}{ll}
            \Theta_\infty(1,\chi)                         &   \ \mbox{if } \chi \mbox{ is of
                                                                                    type } 1 \,, \\
                \ & \\                                                                    
                \displaystyle{\bigg(\frac{\Theta_\infty(\gamma^{-1},\chi)}{1-\gamma^{-1}}\bigg)_{\bigr| \gamma=1} }   &   \ \mbox{if } \chi \mbox{ is of
                                                                                    type } 2 \,,
        \end{array} \right.
    \end{displaymath}
    be the inverse limit of the elements $\Theta_n^{\sharp}(1,\chi)$ appearing in 
    Theorem \ref{Teorema fitting Class Gr per caratteri di tipo 1 o 2}; if $\Theta_\infty^{\sharp}(1,\chi)\neq 0$,
    then $C_\infty(\chi)$ is $\Lambda$-torsion.
\end{teo}

\begin{proof}
   By Proposition \ref{PropNormSurjImmerInject} part (a) and Proposition \ref{PropKernelNormClGr} we
    have
    \begin{displaymath}
        C_m(\chi)=	N^n_m \left(C_n(\chi)\right) \simeq C_n(\chi) / \Ker( N^n_m) = C_n(\chi) /
        I_{\Gamma^n_m}C_n(\chi)\,.
    \end{displaymath}
    Recall that $\mathfrak{I}_m=\displaystyle {\lim_\leftarrow I_{\Gamma_m^n} }$ is the augumentation ideal of
    $\Gal(\mathcal{F} / F_m)$ in $\Lambda$. The previous equality holds for every $n >m$, so
    \begin{displaymath}
        C_m(\chi) \simeq C_\infty(\chi) / \mathfrak{I}_m C_\infty(\chi) \,.
    \end{displaymath}
    The module on the left is a finitely generated torsion $\Lambda / \mathfrak{I}_m=W[\Gamma_m]$-module.
    By the generalized Nakayama Lemma of \cite{BH} we obtain that $C_\infty(\chi)$ is a finitely generated $\Lambda$-module. 
    
    \noindent For the second part just note that
	$\Theta_n^{\sharp}(1,\chi)C_n(\chi) = 0$ for every $n$, hence
    $\Theta_\infty^{\sharp}(1,\chi)C_\infty(\chi) = 0$, and $C_\infty(\chi)$ is a torsion
	$\Lambda$-module.
\end{proof}

\noindent We conclude with the Main Conjecture.

\begin{teo}[\emph{Main Conjecture}]\label{TeoFittIdealClGrIwasawaModule}
    Let $\chi\in\widehat{G_0}$ be a character of type $1$ or $2$. Then 
    \begin{displaymath}
        {\rm Fitt}_{\Lambda} \left(C_\infty(\chi) \right) =
            \left( \Theta_\infty^{\sharp}(1,\chi)\right) \,.
    \end{displaymath}
\end{teo}

\begin{proof}
The proof follows the path of the one of Theorem \ref{TeoFittIdealTateIwasawaModule}.
	The equality
	\[
		\left( \Theta_\infty^{\sharp}(1,\chi)\right) = \lim_\leftarrow
		\left( \Theta_n^{\sharp}(1,\chi)\right) = \lim_\leftarrow 
		{\rm Fitt}_{W[\Gamma_n]} \left( C_n(\chi) \right)
	\]
	reduces the statement to
	\[
		{\rm Fitt}_{\Lambda} \left(C_\infty(\chi) \right) = \lim_\leftarrow  
		{\rm Fitt}_{W[\Gamma_n]} \left( C_n(\chi) \right)\,.
	\]
    Recall $C_m(\chi) = C_\infty(\chi) / \mathfrak{I}_m C_\infty(\chi)$, let $t_1, \dots ,t_r$ be $\Lambda$-generators of $C_\infty(\chi)$
    and denote by $N^\infty_m:C_\infty(\chi) \twoheadrightarrow
    C_m(\chi)$ the projection so that $N^\infty_m(t_1), \dots , N^\infty_m(t_r)$ generate $C_m(\chi)$.
    For any $n$ we have an exact sequence
\begin{equation}\label{EqSeqn}
 \xymatrix{ 0 \ar[r] & K_n \ar[r] & W[\Gamma_n]^{\oplus r} \ar[r] &
C_n(\chi) \ar[r] &  0\,,} 
\end{equation}
   where the map on the right is given by $(w_1, \dots , w_r) \mapsto \sum_i w_i N^\infty_n(t_i)$ and $K_n$ is its kernel.
    They fit into the diagram
\[ \xymatrix{ 0 \ar[r] & K_n \ar[rr] \ar[dd]^{k^n_m} & & W[\Gamma_n]^{\oplus r} \ar[rr] \ar[dd]^{\pi^n_m} 
& & C_n(\chi) \ar[r] \ar[dd]^{N^n_m} & 0 \\
\  \\
0 \ar[r] & K_m \ar[rr] & & W[\Gamma_m]^{\oplus r} \ar[rr]  
& & C_m(\chi) \ar[r] & 0 \,.}\]
    The kernel of $\pi^n_m$ is $(I_{\Gamma^n_m}W[\Gamma_n])^{\oplus r}$
    and, by Proposition \ref{PropKernelNormClGr}, $\Ker(N^n_m)=I_{\Gamma^n_m}C_n(\chi)$, 
    so the map between them is surjective. Moreover $\pi^n_m$ is surjective thus, by the snake lemma, 
    $k^n_m$ is surjective as well and the diagram satisfies the Mittag-Leffler condition. Taking the inverse limit
    and comparing it with the analog of \eqref{EqSeqn} for $K_\infty$ we get $\displaystyle{ K_\infty = \lim_\leftarrow K_n}$.
    
\noindent To conclude the proof, one just follows the same technical arguments of 
	the second part of the proof of Theorem \ref{TeoFittIdealTateIwasawaModule}.
\end{proof}

From this Main Conjecture and the interpolation formulas of the previous sections, one can derive a number of relations with 
special values of $L$ or Zeta-functions. For example we have seen in the 
proofs of Theorems \ref{teorema interpolazione v-adica} and \ref{teorema interpolazione v-adica2},
that for non negative integers $i$ and $j$ with $i\equiv j \pmod{q^{d_\nu}-1}$
\[ \zeta_A(-s_j)(1- \mathfrak{p}^{s_j})=Z(1,j)(1- \mathfrak{p}^{s_j})=L_\mathfrak{p}(1,j,\omega^i)=\prod_{\nu\not\in S} 
(1-\omega(\nu)^i\langle \nu \rangle_{\mathfrak{p}}^j)^{-1}=\Psi_{j,i}(\Theta_S(1)) \,, \]
where the extra factors at the ramified primes different from $\mathfrak{p}$ and $\infty$ disappear 
because in our case $S=\{\mathfrak{p},\infty\}$. Taking $\chi$-parts and projecting from $G_S$ to $\Gamma_\infty$ via $\pi^S_\infty$, one gets
\[ \chi(\pi^S_\infty(L_\mathfrak{p}(X,j,\omega^i)(1)) = \Psi_{j,i}(\Theta_\infty(1,\chi))  \]
for characters of type 1. Up to now  there are, to our knowledge, very few nonvanishing results on
special values for function fields (see e.g. \cite[Theorem E]{AT} and the discussion on ``trivial zeroes'' in \cite[Section 8.13]{Goss libro}): 
it would be interesting
to see if these relations can shed some light on the subject for $L_\mathfrak{p}(1,j,\omega^i)$, which, as seen above, 
is basically equivalent to $C_\infty$ being a torsion $\Lambda$-module. 

For $F=\F_q(t)$, in \cite[Section 6]{ABBL} the authors present arithmetic information on some Bernoulli-Goss numbers, i.e. on special
values of the Goss Zeta-function at integers: in our setting one should probably consider the finite $\F_q[t]$-module $H(\Phi/A)$ defined by 
Taelman in \cite{Tael1} (where $\Phi$ is a Drinfeld module over $A$), which plays the role of the ideal class group of a number field. 
For $F=\F_q(t)$ the Bernoulli-Goss numbers are linked to the isotypical components of $H(\Phi/A)$ by \cite[Theorem 1 and Section 10]{Tael2}. 
It would be interesting to study inverse limits of 
$\mathfrak{p}$-parts of  Taelman's modules associated to the fields $F_n$ as objects over the Iwasawa algebra:
it is not clear whether this would lead to special values of our $\mathfrak{p}$-adic $L$-function or of some other (yet to be defined) 
$\mathfrak{p}$-adic analytic function (possibly another incarnation of the Stickelberger series).

\medskip

\noindent {\bf Acknowledgements.} We are grateful to Bruno Angl\'es, Francesc Bars and Ignazio Longhi for several useful conversations, 
suggestions and comments which provided invaluable contributions to the development of this paper. We would like to thank
Fabrizio Andreatta and Marco Seveso for their help and support. We also thank the anonymous referees for 
their remarks which improved the exposition and provided inputs for future research.


\begin{thebibliography}{99}
	
\bibitem{AT} B. Angl\'es - L. Taelman, ``Arithmetic of characteristic $p$ special $L$-values''
(with an appendix by V. Bosser)
{\em Proc. Lond. Math. Soc.} (3) {\bf 110} (2015), 1000--1032.	

\bibitem{ABBL} B. Angl\'es - A. Bandini - F. Bars - I. Longhi, ``Iwasawa Main Conjecture for the Carlitz cyclotomic 
extension and applications'', {\em Math. Ann.} {\bf 376}, Issue 1-2 (2020), 475--523.

\bibitem{BH} P.N. Balister - S. Howson, ``Note on Nakayama's lemma for compact $\Lambda$-modules'',
{\em Asian J. Math.} {\bf 1} (1997), no. 2, 224--229.

\bibitem{BBC} A. Bandini - F. Bars - E. Coscelli, ``Fitting ideals of class groups in Carlitz-Hayes cyclotomic extensions'',
to appear in {\em J. Number Theory}.

\bibitem{BBL} A. Bandini - F. Bars - I. Longhi, ``Characteristic ideals and Iwasawa theory'',
{\em New York J. Math} {\bf 20} (2014), 759--778.

\bibitem{BBL2} A. Bandini - F. Bars - I. Longhi, ``Characteristic ideals and Selmer groups'',
{\em J. Number Theory} {\bf 157} (2015), 530--546.

\bibitem{BV} A. Bandini - M. Valentino, ``Control theorems for $\ell$-adic Lie extensions of global function fields'',
{\em Ann. Sc. Norm. Super. Pisa Cl. Sci.}  {\bf XIV} (2015), no. 4, 1065--1092. 

\bibitem{BV2} A. Bandini - M. Valentino, ``Euler characteristic and Akashi series for Selmer groups 
over global function fields'',
{\em J. Number Theory} {\bf 193} (2018), 213--234. 

\bibitem{Bur} D. Burns, ``Congruences between derivatives of geometric L-functions.'' 
With an appendix by Burns, K.F. Lai and K.-S. Tan, {\em Invent. Math.} 
{\bf 184} (2011), no. 2, 221--256.

\bibitem{Carlitz} L. Carlitz, ``On certain functions connected with polynomials in a Galois field'', 
\emph{Duke Math. J.} (1935), 137--168.

\bibitem{Cornacchia Greither} P. Cornacchia - C. Greither, 
``Fitting ideals of class group of real fields with prime power conductor'', 
\emph{J. Number Theory} \textbf{73} (1998), 459--471.

\bibitem{CoscPhD} E. Coscelli, ``Stickelberger series and Iwasawa Main Conjecture for function fields'', PhD Th., 
University of Milan (2018), avaliable at \text{https://air.unimi.it/handle/2434/561439\#.XUWhDHvONPY} .

\bibitem{Crew} R. Crew, ``$L$-functions of $p$-adic characters and geometric Iwasawa theory'',
{\em Invent. Math.} {\bf 88} (1987), no. 2, 395--403.

\bibitem{Goss libro} D. Goss, \emph{Basic Structures of Function Field Arithmetic,} (Springer-Verlag, 1996).

\bibitem{Goss Zeta} D. Goss, ``$v$-adic Zeta-functions, $L$-series and Measures for Function Fields'', 
\emph{Invent. Math.} \textbf{55} (1979), 107--116.

\bibitem{Greither Kurihara} C. Greither - M. Kurihara,
``Stickelberger elements, Fitting ideals of class groups of CM-fields and dualisation'', 
\emph{Math. Z.} \textbf{260} (2008), no. 4, 905--930.

\bibitem{Greither Popescu Galois} C. Greither - C.D. Popescu,
``The Galois module structure of $\ell$-adic realizations of Picard $1$-motives and applications'', 
\emph{Int. Math. Res. Not.} (2012), no. 5, 986--1036.

\bibitem{Greither Popescu Fitting} C. Greither - C.D. Popescu, 
``Fitting ideals of $\ell$-adic realizations of Picard $1$-motives and class groups of global function fields'', 
\emph{J. Reine Angew. Math.}  \textbf{675} (2013), 223--247.

\bibitem{Hayes Drinfeld} D.R. Hayes, ``A brief introduction do Drinfeld Modules'' in ``The Arithmetic of function fields'' 
(\emph{Columbus,} OH, 1991) Ohio State Univ. Math. Res. Inst. Publ. \textbf{2}, 1--32.

\bibitem{KatoICM} K. Kato, ``Iwasawa theory and generalizations'', ICM Vol I, Eur. Math. Soc.
Z\"urich, Vol. 12, 335--357 (2007).


\bibitem{LLTT} K.-F. Lai - I. Longhi - K.-S. Tan - F. Trihan, ``The Iwasawa main conjecture for constant ordinary abelian varieties 
over function fields'', 
{\em Proc. Lond. Math. Soc.} {\bf 112} (2016), , no. 6, 1040--1058.

\bibitem{Mazur Wiles} B. Mazur - A. Wiles,
``Class fields of abelian extensions of $\Q$'', 
\emph{Invent. Math.} \textbf{76} (1984), 179--330.

\bibitem{Mumford} D. Mumford, {\em Abelian Varieties}, 
Tata Inst. of Fundamental Research, Bombay, 1970.

\bibitem{Northcott} D.G. Northcott, {\em Finite free resolutions},
Cambridge University Press, Cambridge Tracts in Mathematics, No. \textbf{71}, Cambridge, 1976.

\bibitem{Rosen} M. Rosen, {\em Number theory in function fields},
{\bf GTM 210}, Springer-Verlag, New York, 2002.


\bibitem{Shu Kummer} L. Shu, ``Kummer's criterion over global function fields'', 
{\em J. Number Theory} {\bf 49} (1994), 319--359.

\bibitem{Tael1} L. Taelman, ``Special $L$-values of Drinfeld modules'',
{\em Annals of Math.} {\bf 175} (2012), 369--391.

\bibitem{Tael2} L. Taelman, ``A Herbrand-Ribet theorem for function fields'',
{\em Invent. Math.} {\bf 188} (2012), 253--275.

\bibitem{Tan} K.-S. Tan, ``A generalized Mazur’s theorem and its applications'', 
{\em Trans. Amer. Math. Soc.} {\bf 362} (2010), 4433--4450.

\bibitem{Tan2} K.-S. Tan, ``Selmer groups over $\Z_p^d$-extensions'', 
{\em Math. Ann.} {\bf 359} (2014), 1025--1075.


\bibitem{TateStark} J. Tate,
{\em Les conjectures de Stark sur les Fonctiones $L$ d'Artin en $s=0$}, 
Progress in Mathematics \textbf{47}, Birkh\"auser, (1984).

\bibitem{Thakur} D.S. Thakur, {\em Function Field Arithmetic},
World Scientific Publishing Co., Inc., River Edge, NJ, 2004.


\end{thebibliography}
\end{document}